\documentclass[11pt,a4paper]{article}

\usepackage{psfrag}
\usepackage{amsmath,amssymb,cite}
\usepackage{latexsym}
\usepackage{graphicx}
\usepackage{lscape}
\usepackage{afterpage}

\makeatother\makeatother\usepackage[color=green!40]{todonotes}

\usepackage{hyperref}
\hypersetup{
    colorlinks=true,
    linkcolor=blue,
    filecolor=magenta,      
    urlcolor=cyan,
    citecolor=red,
}
\makeatletter
\let\atop\@@atop
\makeatother

\DeclareMathOperator*{\argmin}{argmin}

\newcommand{\varphib}{\overline{\varphi}}
\newcommand{\xb}{\overline{x}}
\newcommand{\ub}{\overline{u}}
\newcommand{\tb}{\overline{t}}

\newcommand{\ds}{\displaystyle}
\newcommand{\nexto}{\kern -0.54em}
\newcommand{\dR}{{\rm {I\ \nexto R}}}

\newcommand{\dZ}{{\cal Z \kern -0.7em Z}}
\newcommand{\dC}{{\rm\hbox{C \kern-0.8em\raise0.2ex\hbox{\vrule
height5.4pt width0.7pt}}}}
\newcommand{\dQ}{{\rm\hbox{Q \kern-0.85em\raise0.25ex\hbox{\vrule
height5.4pt width0.7pt}}}}
\newcommand{\proofbox}{\hspace{\fill}{$\Box$}}

\newtheorem{theorem}{Theorem}

\newtheorem{fact}{Fact}
\newtheorem{remark}{Remark}

\newtheorem{algorithm}{Algorithm}

\begin{document}

\author{
C. Yal{\c c}{\i}n Kaya\footnote{Mathematics, UniSA STEM, University of South Australia, Mawson Lakes, S.A. 5095, Australia. E-mail: yalcin.kaya@unisa.edu.au\,.}
\and
Helmut Maurer\footnote{Institut f\"ur Numerische und Angewandte Mathematik, Westf\"alische Wilhelms-Universit\"at M\"unster, M\"unster, Germany. E-mail: helmut.maurer@uni-muenster.de\,.}
}

\title{\vspace{-10mm}\bf Optimization Over the Pareto Front of Nonconvex Multi-objective Optimal Control Problems}

\maketitle

\thispagestyle{empty}

\begin{abstract} {\noindent\sf 
Simultaneous optimization of multiple objective functions results in a set of trade-off, or Pareto, solutions. Choosing a, in some sense, best solution in this set is in general a challenging task: In the case of three or more objectives the Pareto front is usually difficult to view, if not impossible, and even in the case of just two objectives constructing the whole Pareto front so as to visually inspect it might be very costly.  Therefore, optimization over the Pareto (or efficient) set has been an active area of research.  Although there is a wealth of literature involving finite dimensional optimization problems in this area, there is a lack of problem formulation and numerical methods for optimal control problems, except for the convex case. In this paper, we formulate the problem of optimizing over the Pareto front of nonconvex constrained and time-delayed optimal control problems as a bi-level optimization problem. Motivated by existing solution differentiability results, we propose an algorithm incorporating (i) the Chebyshev scalarization, (ii) a concept of the essential interval of weights, and (iii) the simple but effective bisection method, for optimal control problems with two objectives.  We illustrate the working of the algorithm on two example problems involving an electric circuit and treatment of tuberculosis and discuss future lines of research for new computational methods.}
\end{abstract}

\begin{verse}
{\em Key words}\/: {\sf Multi-objective optimization, Optimal control, Optimization over Pareto front, Optimization over efficient set, Numerical methods, Rayleigh problem, Tuberculosis, Time-delay problems.}
\end{verse}

%\centerline{\bf Submitted to }
\pagestyle{myheadings}
%\markboth{}{Submitted to {\sf }, }
\markboth{}{\sf\scriptsize Optimization Over the Pareto Front of
  Multi-objective Optimal Control Problems\ \ by C. Y. Kaya and
  H. Maurer}

\section{Introduction}

We continue our study of optimal control problems where one wishes to minimize {\em simultaneously} a number of conflicting objective functionals. These problems are referred to as multi-objective optimal control problems and can be expressed in the following concise form:
\[
\mbox{(P)}\hspace*{10mm} \min_{(x,u,t_f)\in X}\ (\varphi_1(x(t_f),t_f),\ldots,\varphi_r(x(t_f),t_f))\,.
\]
The constraint or the feasible set $X$ in Problem~(P) involves a system of differential equations (DEs) in the state and control variables $x(\cdot)$ and $u(\cdot)$, respectively, over a time horizon $[0,t_f]$.  The set $X$ also typically involves point and path equality and inequality constraints.  The DEs and constraints in $X$ might even include time delays in the variables $x(\cdot)$ and $u(\cdot)$.  It is worth noting that although each of the objective functionals $\varphi_i(x(t_f),t_f)$, $i = 1,\ldots,r$, in~(P) above constitutes the so-called {\em Mayer form}, other forms (Bolza and Lagrange) can easily be converted into this form conveniently.  Therefore, the general model in (P) caters for a wide range of conflicting objectives; for instance, minimization of the energy, the terminal time, the deviations from a reference state trajectory, or the uncertainty in measurements, to name just a few.

Broadly speaking, the {\em simultaneous} or {\em Pareto minimization} in Problem~(P) is the process of finding a compromise solution, referred to as a {\em Pareto minimum}, where the value of some cost cannot be improved (i.e., reduced) further, without making the value of some other cost worse (i.e., higher).  One typical example is the case when one wants to minimize simultaneously the fuel expenditure of an airplane travelling from one given city to another and the time the airplane takes for this travel: A shorter travel time often requires a higher fuel consumption.   The set of all such compromise or trade-off solutions form the {\em Pareto set} in the optimization space, or the {\em Pareto front} in the value space.  Pareto set and Pareto front are also commonly referred to as the {\em efficient set} and the {\em efficient front}, respectively\footnote{These and other definitions will be given in more precise terms in Section~\ref{sec:prelim}.}.

The authors of this paper have studied in~\cite{KayMau2014} the problem of constructing the Pareto front of Problem~(P) involving ODEs and constraints of general form.  They discussed and demonstrated that for the nonconvex optimal control problems like the one in Problem~(P), it is better to use the so-called weighted Chebyshev-norm scalarization (or just Chebyshev scalarization) to guarantee that the whole Pareto front can be constructed, instead of using the traditional weighted-sum scalarization, i.e., a convex combination of the objective functionals.  They discretized the scalarized problem directly and utilized large-scale optimization software (the AMPL--Ipopt suite~\cite{AMPL, WacBie2006}) to find the Pareto fronts of two  constrained optimal control problems as examples, one involving tumour anti-angiogenesis and the other a fed-batch bioreactor, by means of what they called a {\em scalarize--discretize--then--optimize} approach.  This approach is in contrast with the other existing {\em discretize--scalarize--then--optimize} approach (see e.g. \cite{LogErdImp2010, LogHouDieImp2010, LogValHouDieImp2012, ObeRinFel2012}) which scalarizes the discretized problem rather than the original (continuous-time) problem.

An additional benefit of the Chebyshev scalarization is also reported and illustrated in~\cite{KayMau2014}:  One can compute the whole Pareto front by using only those weights of the objective functionals within what they name as the {\em essential subinterval} of weights, instead of the whole interval.  Having to compute fewer Pareto solutions over a smaller number of grid points in a subinterval is obviously a computational advantage.  For further details and an extensive list of references on multi-objective optimal control the reader is referred to~\cite{KayMau2014}.  Other relevant studies on the topic in more recent years have appeared in~\cite{DesHas2019, Chorobura2021}.

Apart from certain trivial or special cases, the Pareto front consists of infinitely many solutions to choose from.  When a discrete approximation of the front is found the number of solutions to choose from is still relatively large since the approximate front is required to be accurate enough. Making a decision as to which Pareto solution in the front is the most suitable (to the needs of a practitioner) is often very hard for the following reasons.
\begin{itemize}
\item  In the case of three or more objectives, the Pareto front might be difficult (if not impossible) to  view and to carry out a visual inspection (or ``eyeballing'') for a decision.
\item  Even with two objectives, a visual inspection alone may not be enough to choose a desirable solution.
\item  Constructing the whole Pareto front might just be too costly a thing to do numerically.
\end{itemize}
Motivated by these drawbacks, minimization of an additional (single) objective function over the Pareto front has been of great interest to many researchers over the past decades---see, for example, \cite{Benson1984, Benson1992, Bolintineanu1993a, Bolintineanu1993b, Dauer1991, DauFos1995, HorTho1999, HorThoYamZen2007, LiuEhr2018, Philip1972, Yamamoto2002}.  Despite this rich collection of works, to the knowledge of the authors, it was not before the reference~\cite{BonKay2010} that optimization over the Pareto front was studied and a numerical method proposed for {\em convex} multi-objective optimal control problems.  In the current paper, we extend the works in~\cite{BonKay2010,KayMau2014} to {\em nonconvex} multi-objective optimal control problems and propose a numerical method for carrying out optimization over the Pareto front.

We set the optimal control problem as a bi-level optimization problem as in~\cite{BonKay2010}: One has to minimize a master objective functional subject to the minimization of a scalarization of Problem~(P).  The lower level problem uses the Chebyshev scalarization as in~\cite{KayMau2014}, as opposed to the weighted-sum scalarization in~\cite{BonKay2010}.  The problems we consider is in much more general form in this paper: We consider nonconvex instead of convex problems compared to~\cite{BonKay2010} and we consider problems with time-delay instead of those without time delay compared to~\cite{KayMau2014}.  Just to re-iterate, \cite{KayMau2014} only proposes a technique to construct the Pareto front, otherwise it does not carry out optimization over the Pareto front.

As the optimization technique over the Pareto front, we propose the simplest possible technique, namely the bisection method, over the set of weights for the bi-objective problem, which are the parameters of the lower level optimal control problem.  Even in this simplest case, it is necessary to obtain derivatives with respect to the weight, for which we employ difference approximations.  However, is it guaranteed that these derivatives exist?  This question is answered by~\cite{Malanowski-Maurer-96,Malanowski-Maurer-98,MauPes1994,MauPes1995} which studied the differentiability of a solution of a parametric optimal control problem with respect to the parameters.  We add a discussion concerning these studies in the paper.

The main algorithm first finds the essential interval of weights over which the first step of the bisection method is taken to find a new subinterval.  Then the subsequent steps of the bisection method are carried out until the stopping criterion is met.

The algorithm is illustrated on two challenging numerical examples: the Rayleigh problem, which comes from an electric circuit, and a compartmental optimal control model for tuberculosis.  In the first problem there are constraints on the control variables, and the second problem not only has constraints on the two control variables but also time delays on both the control and state variables.

The paper is organized as follows.  In Section~2, we introduce the multi-objective optimal control problem, discuss scalarization, introduce the problem of optimization over the Pareto front, and elaborate on solution differentiability.  In Section~3, we first define and explain the essential interval of weights, and then introduce the bisection method for our problem and provide the detailed algorithm.  In Section~4, we illustrate the algorithm on two example optimal control problems.  Finally, in Section~5, we provide concluding remarks.

\section{Problem Statement and Preliminaries}
\label{sec:prelim}

\subsection{Multi-objective optimal control problem}
We consider the following general multi-objective optimal control problem (similar to that in~\cite{KayMau2014} but made look slightly more general here) to underlie our study on  minimization over its Pareto front.  The ensuing notation and definitions can also be found in~\cite{KayMau2014} but given here for completeness as well as convenience.
\[
\mbox{(OCP) }\left\{\begin{array}{rl}
\ds\min & \ \ \ds
(\varphi_1(x(t_f),t_f),\ldots,\varphi_r(x(t_f),t_f)) \\[3mm]
\mbox{subject to} & \ \ \dot x(t) = f(x(t),u(t),t)\,, \mbox{\ \ for a.e. } t\in[0,t_f]\,, \\[2mm]
& \ \ \theta(x(0),x(t_f),t_f) = 0\,, \\[2mm]
& \ \ \widetilde\theta(x(0),x(t_f),t_f) \le 0\,, \\[2mm]
& \ \ C(x(t),u(t),t) \le 0\,, \mbox{\ \ for a.e. } t\in[0,t_f]\,, \\[2mm]
& \ \ S(x(t),t) \le 0\,, \mbox{\ \ for all } t\in[0,t_f]\,,
\end{array} \right.
\]
where $r\in\{2,3,4,\ldots\}$ is fixed, the state variable $x\in W^{1,\infty}(0,t_f;\dR^n)$, $\dot{x} := dx/dt$, and the control variable $u\in L^{\infty}(0,t_f;\dR^m)$, with $x(t) := (x_1(t),\ldots,x_n(t))\in\dR^n$ and $u(t) := (u_1(t),\ldots,u_m(t))\in\dR^m$.  The functions $\varphi_i:\dR^n\times\dR_+\rightarrow\dR$,
$f:\dR^n\times\dR^m\times\dR_+\rightarrow\dR^n$,
$\theta:\dR^n\times\dR^n\times\dR_+\rightarrow\dR^{p_1}$,
$\widetilde\theta:\dR^n\times\dR^n\times\dR_+\rightarrow\dR^{p_2}$,
$C:\dR^n\times\dR^m\times\dR_+\rightarrow\dR^{p_3}$, and
$S:\dR^n\times\dR_+\rightarrow\dR^{p_4}$, 
are continuous in their arguments.  In this problem, $t_f$ is either fixed or free.  Here, $L^{\infty}(0,t_f;\dR^m)$ corresponds to the space of essentially bounded, measurable functions equipped with the essential supremum norm. Furthermore, $W^{1,\infty}(0,t_f;\mathbb{R}^n)$ is the Sobolev space consisting of functions $x:[0,t_f]\rightarrow \mathbb{R}^n$ whose first derivatives lie in $L^{\infty}$.

Assume that $\varphi_i(x(t_f),t_f)\ge 0$, for all $i=1,\ldots,r$.  Note that this assumption can easily be met by adding a large enough positive number to each objective functional.

Note that Problem~(OCP) is in general a nonsmooth problem, because it does not require differentiability of the objective functionals or the constraints.  Moreover, although we have stated Problem~(OCP) in very broad terms, it can further be generalized, for example by adding multi-point constraints, partial differential equations, time delays, etc.  In other words, although Problem~(OCP) is already in a more general form than what one usually encounters in applications, it can be further made look more general.

Of the possible extensions mentioned above, {\em time delays} in the state and control variables, for instance, can be incorporated into Problem~(OCP) by replacing the ODEs in Problem~(OCP) with
\begin{subequations}
\begin{eqnarray}  \label{delayed_ODE}
&& \dot x(t) = f(x(t),x(t-d_x),u(t),u(t-d_u),t)\,, \mbox{\ \ for a.e. } t\in[0,t_f]\,, \label{delayed_ODEa} \\[1mm]
&& x(t) = x_0(t)\,, \mbox{\ \ for all } t\in[-d_x,0)\,, \label{delayed_ODEb} \\[1mm]
&& u(t) = u_0(t)\,, \mbox{\ \ for all } t\in[-d_u,0)\,, \label{delayed_ODEc}
\end{eqnarray}
\end{subequations}
where $d_x,d_u > 0$ are the time delays in the state and control variables, respectively.

For technical convenience, let $t_f\le t_f^{\max}$, where $t_f^{\max}>0$ is some constant.  Next, we define the {\em feasible set}, $X\subset W^{1,\infty}(0,t_f;\dR^n)\times L^\infty(0,t_f;\dR^m)\times\dR_+$, such that
\begin{eqnarray*}
X &:=&\{(x,u,t_f)\,:\,\dot x(t) = f(x(t),x(t-d_x),u(t),u(t-d_u),t)\,, \mbox{\ \ for a.e. } t\in[0,t_f]\,; \\[1mm]
&& \quad x(t) = x_0(t)\,,\mbox{\ \ for all } t\in[-d_x,0];\ \ u(t) = u_0(t)\,,\mbox{\ \ for all } t\in[-d_u,0)\,; \\[1mm]
&& \quad \theta(x(0),x(t_f),t_f) = 0\,;\ \widetilde\theta(x(0),x(t_f),t_f) \le 0\,; \\[1mm]
&&\quad C(x(t),u(t),t) \le 0\,,\mbox{ for a.e. } t\in[0,t_f];\ S(x(t),t) \le 0,\mbox{ for all } t\in[0,t_f]\}\,.
\end{eqnarray*}
Note that, for the case of time delays in the state and control variables, we have included Equations~\eqref{delayed_ODEa}--\eqref{delayed_ODEc} instead of the ODEs\ \ $\dot{x}(t) = f(x(t),u(t),t)$ in the set $X$.

Define the vector of objective functionals, $\varphi(x(t_f),t_f) := (\varphi_1(x(t_f),t_f),\ldots,\varphi_r(x(t_f),t_f))$.  The triplet $(x^*,u^*,t_f^*)\in X$ is said to be a {\em Pareto minimum} if there exists no $(x,u,t_f)\in X$ such that $\varphi(x(t_f),t_f) \neq \varphi(x^*(t_f^*),t_f^*)$ and
\[
\varphi_i(x(t_f),t_f)) \le \varphi_i(x^*(t_f^*),t_f^*)\,,\quad\mbox{for all } i=1,\ldots,r\,.
\]
On the other hand, $(x^*,u^*,t_f^*)\in X$ is said to be a {\em weak Pareto minimum} if there exists no $(x,u,t_f)\in X$ such that
\[
\varphi_i(x(t_f),t_f)) < \varphi_i(x^*(t_f^*),t_f^*)\,,\quad\mbox{for all } i=1,\ldots,r\,.
\]
The set of all the Pareto and weak Pareto minima is said to be the {\em Pareto set}.  On the other hand, the set of all vectors of objective functional values at the Pareto and weak Pareto minima is said to be the {\em Pareto front} (or the {\em efficient set}) of Problem~(OCP) in the $r$-dimensional\linebreak {\em objective value}, or {\em outcome}, {\em space}.  Note that the coordinates of a point in the Pareto front are simply $\varphi_i(x^*(t_f^*),t_f^*)$, $i=1,\ldots,r$.  Obviously, when $r=2$ the Pareto front is in general a curve; and when $r=3$ the Pareto front is in general a surface.

\subsection{Scalarization}
\label{scalarization}

In~\cite{KayMau2014}, to compute a solution of Problem (OCP), the following single-objective problem~(P$_w$), i.e., {\em scalarization}, was employed.
\[
\mbox{(P$_w$)}\qquad \min_{(x,u,t_f)\in X}\ \max\{w_1\,\varphi_1(x(t_f),t_f),\ldots,w_r\,\varphi_r(x(t_f),t_f)\}\,,
\]
where $w_i$, $i=1,\ldots,r$, are referred to as {\em weights}, with the vector of weights $w$ defined as $w := (w_1,\ldots,w_r)\in\dR^r$, such that $\sum_{i=1}^r w_i = 1$.  Problem~(P$_w$) is referred to as the {\em weighted Chebyshev problem} (or {\em Chebyshev scalarization}) because of the weighted Chebyshev norm, $\max_i |w_i\,\varphi_i(x(t_f),t_f)|=\max_i w_i\,\varphi_i(x(t_f),t_f)$, appearing in the objective.  This type of scalarization is typically used for nonconvex multi-objective finite-dimensional optimization problems, as opposed to the weighted sum scalarization which is effective for convex problems but not the nonconvex ones---see, for example, \cite{Miettinen1999}.

Define the set of weights
\[
Y := \left\{w\in\dR^r\ |\ \sum_{i=1}^r w_i = 1\right\}\,.
\]
The following theorem was originally presented in~\cite[Theorem~1]{KayMau2014} for the case when there was no delay in the state and control variables.  It still holds with the set $X$ modified with the delayed state equations.

\begin{theorem} [Bijection between sets of weights and Pareto minima~\cite{KayMau2014}]  \label{scalar_theorem}
The triplet $(x^*,u^*,t_f^*)$ is a weak Pareto minimum of~(OCP) if, and only if, $(x^*,u^*,t_f^*)$ is a solution of {\rm (P$_w$)} for some $w_1,\ldots,w_r>0$.
\end{theorem}

\begin{remark}  \label{rem:surjection} \rm
Suppose that $Z\subset X$ denotes the Pareto set, namely the set of all Pareto minima of (OCP).  Then Theorem~\ref{scalar_theorem} establishes that there is a bijection between the set of weights $Y$ and the Pareto set $Z$.  This implies that by solving (P$_w$) for all $w\in Y$, one can obtain the whole Pareto set $Z$ and in turn get the Pareto front.  With numerical computations on the other hand, one would of course carry out some discretization of the weight space $Y$ and typically get a discrete approximation of the Pareto front.  The bijection between $Y$ and $Z$ will also help us devise our algorithm for optimization over the Pareto front.
\proofbox
\end{remark}

An {\em ideal cost} $\varphi^*_i$, $i=1,\ldots,r$, associated with Problem~(P$_w$) is the optimal value of the optimal control problem,
\begin{equation}  \label{single_obj}
\min_{(x,u,t_f)\in X}\ \varphi_i(x(t_f),t_f)\,.
\end{equation}
Let $(\xb,\ub,\tb_f)$ be a minimizer of the single-objective problem in~\eqref{single_obj}.  Then $\varphi^*_i := \varphi_i(\xb(\tb_f),\tb_f)$ and we also define $\varphib_j :=
\varphi_j(\xb(\tb_f),\tb_f)$, for $j\neq i$ and $j=1,\ldots,r$.

In the case when $\varphi^*_i$ is negative, one can simply add a large enough positive number to the $i$th objective, to make the objective positive.  In general, it is useful to add a positive number to each
objective in order to obtain an even spread of the Pareto points approximating the Pareto front -- see for example \cite{DutKay2011} for further discussion and geometric illustration.  To serve this purpose, it is common practice to define the so-called {\em utopian objective values}.

A {\em utopian objective vector} associated with Problem~(OCP) is given as $\beta^* := (\beta_1^*,\ldots,\beta_r^*)$, with $\beta^*_i := \varphi_i^* - \eta_i$ and $\eta_i > 0$ for all $i = 1, \ldots, r$.  Problem~(P$_w$) can then be equivalently written as
\[
\min_{(x,u,t_f)\in X}\ \max\{w_1\,(\varphi_1(x(t_f),t_f) - \beta_1^*),\ldots,w_r\,(\varphi_r(x(t_f),t_f) - \beta_r^*)\}\,.
\]

In the case when the objective functionals and the constraints in Problem~(OCP) are differentiable in their arguments, it is worth reformulating Problem~(P$_w$) using a standard technique from mathematical programming in the following (smooth) form.
\[
\mbox{(OCP$_w$) }\left\{\begin{array}{rl}
\ds\min_{{\alpha\ge0}\atop{(x,u,t_f)\in X}} & \ \ \alpha \\[4mm]
\mbox{subject to} & \ \ w_1\,(\varphi_1(x(t_f),t_f) - \beta_1^*) \le\alpha\,, \\
			  & \hspace*{24mm}\vdots \\
			  & \ \ w_r\,(\varphi_r(x(t_f),t_f) - \beta_r^*) \le\alpha\,.
\end{array} \right.
\]
Problem (OCP$_w$) is referred to as {\em goal attainment method} \cite{Miettinen1999}, as well as {\em Pascoletti-Serafini scalarization} \cite{Eichfelder2008}.  We will solve Problem~(OCP$_w$) in an algorithm we present in the next section, for the two examples we want to study.

We re-iterate that the ``popular'' {\em weighted-sum} scalarization, given below, fails to generate the ``nonconvex parts'' of a Pareto front.
\[
\mbox{(P$_{ws}$)}\qquad \min_{(x,u,t_f)\in X}\ \sum_{i=1}^{r} w_i\,\varphi_i(x(t_f),t_f)\,.
\]
This deficiency is illustrated with a multi-objective optimal control problem, for example, in the fed-batch bioreactor problem in~\cite{KayMau2014}.

\subsection{Optimization over the Pareto front}
\label{sec:OPF}

The main task in this paper is to devise a numerical algorithm for solving the problem of decision making as to which Pareto point should be chosen.  This obviously depends on the criterion a decision maker uses in making his/her choice.  As pointed in Remark~\ref{rem:surjection}, the whole Pareto front can be parameterized in terms of the vector of weights $w$.  Therefore, Problem~(P$_w$), or equivalently (OCP$_w$), can be regarded as a parametric optimal control problem, and it also makes sense to express the decision maker's objective as the minimization of a function of $w$.  

Before going ahead with the statement of this problem, we re-write the variables of the optimal control problem, with a slight abuse of notation, as $x^w(t) := x(t,w)$, $u^w(t) := u(t,w)$, and $t_f^w := t_f(w)$ to emphasize their dependence on the vector of weights $w$.  
%The terminal time $t_f$, when free, also depends of course on $w$; however, we do not show that explicitly here for the sake of clarity in the appearance of the expressions.  
%When $t_f$ is fixed, as it happens to be the case with many optimal control problems, showing dependence on $w$ does not make so much sense. % anyway.

We call the decision maker's objective function the {\em master objective function}, expressed by $\varphi_0(x^w, u^w, t_f^w)$.  With the weight vector $w$ of the scalarization treated now as a variable, the problem of optimization over the Pareto front reduces to the problem of finding an optimal weight~$w^*$.  Then the corresponding Pareto minimum is a solution of Problem~(OCP$_{w^*}$).  

The problem of optimizing a master objective function over the Pareto front of~(OCP) with $r\ge2$ objectives is nothing but a bilevel programming problem and can be written as
\[
\mbox{(OPF) }\left\{\begin{array}{rll}
\ds\min_{w\in Y} & \ \ \varphi_0(x^w, u^w, t_f^w) & \\[2mm]
 \ \ \mbox{subject to}  & \ \ \ds\min_{{\alpha\ge0}\atop{(x,u,t_f)\in X}} & \hspace*{-7mm} \alpha \\[4mm]
 & \ \ \mbox{subject to} & \hspace*{-7mm} w_1\,(\varphi_1(x(t_f,w),t_f) - \beta_1^*)
\le\alpha\,, \\
 & \ \ & \hspace*{15mm}\vdots \\
 & \ \ & \hspace*{-7mm} w_r\,(\varphi_r(x(t_f,w),t_f) - \beta_r^*) \le\alpha\,.
\end{array} \right.
\]

\begin{remark}  \rm
The lower-level problem in (OPF) for some given $w$ is simply Problem~(OCP$_w$).  A solution of~(OCP$_w$) is nothing but a point in the Pareto set $Z$ of (OCP) and is described by the triplet $Z_w := (x^*(t,w), u^*(t,w), t_f^*(w))$.  Then the (whole) Pareto set can be expressed as $Z = \cup_{w\in Y} Z_w$. Now Problem~(OPF) can equivalently be written as
\[
\left\{\begin{array}{rll}
\ds\min_{w\in Y} & \ \ \varphi_0(x^w, u^w, t_f^w) & \\[2mm]
 \ \ \mbox{subject to}  & \ \ (x^w, u^w, t_f^w)\in Z_w\,.
\end{array} \right.
\]
We note that the optimization variable of the upper-level problem is the ``unknown'' parameter $w$.  If the solution $(x^*(t,w), u^*(t,w), t_f^*(w))$ of Problem~(OCP$_w$) is differentiable in the parameter $w$, then powerful differentiable optimization techniques can be employed in solving Problem~(OPF) (or in a more concise form the above problem).  This is what was done in~\cite{BonKay2010} for convex multi-objective optimal control problems. In this paper, we are extending the work in~\cite{BonKay2010} to the nonconvex setting by also incorporating the Chebyshev scalarization and the concept of essential interval of weights given in~\cite{KayMau2014}.
\proofbox
\end{remark}

\subsection{Solution differentiability}
\label{sec:soldiff}

% Maurer and Pesch study in~\cite{MauPes1995} the differentiability
% of a local solution of the following class of parametric optimal
% control problems.
We briefly review results on solution differentiability or $C^1$-sensitivity of solutions to
 the following parametric optimal control problems depending on a parameter 
% $p \in \mathbb{R}^s$
$p \in P$, where $P$ is a Banach space:
\[
\mbox{(OCP(p))}
\left\{\begin{array}{rl}
\min_{\,x,u,p} &  g(x(t_f),t_f,p) \\[2mm]
\mbox{subject to} &  \dot x(t) = \widetilde{f}(x(t),u(t),p)\,, \mbox{\ \ for a.e. } t\in[0,t_f]\,, \\[1mm]
&  \psi(x(0),x(t_f),t_f,p)  = 0 \,,\quad \tilde{\psi}(x(0),x(t_f),t_f,p) \leq 0 \,,  \\[1mm]
& \tilde{C}(x(t),u(t),p) \le 0\,, \mbox{\ \ for a.e. } \; t\in[0,t_f]\,, \\[1mm]
& \tilde{S}(x(t),p) \le 0\,, \mbox{\ \ for a.e. } \;  t\in[0,t_f]\,.
\end{array} \right.
\]
We note that problem (OCP$_w$) is a special case of the parametric problem (OCP(p)) by simply taking the parameter 
as the weight, $p=w$, which then appears only in the terminal inequality constraints.
The problem (OCP($p_0$)) corresponding to a reference  parameter $p_0$ is considered as the \textit{nominal} or
\textit{unperturbed} problem. It is assumed that a local solution $(x_0,u_0)$ of the reference solution
exists. 
Let $p$ be a parameter in a neighbourhood of the nominal parameter $p_0$ and denote the
solution to (OCP($p$)) by $(x(t,p),u(t,p))$. 
Dontchev and Hager \cite{Dontchev-Hager-93} gave conditions under which the mapping 
$p\mapsto (x(\cdot,p), u(\cdot,p))$ is Lipschitz. 
Malanowski and Maurer \cite{Malanowski-Maurer-96,Malanowski-Maurer-98} and Maurer and Pesch \cite{MauPes1994,MauPes1995} investigated the \textit{solution differentiability} or $C^1$-sensitivity of the
 optimal solution. The authors derived conditions such that
 % %  showed that under strong second-order sufficient conditions for the nominal solution
 an optimal solution $(x(\cdot,p),u(\cdot,p))$ of the perturbed control problem OCP($p$) exists for all parameters $p$ in a neighborhood of $p_0$ and, moreover, the solution $(x(t,p),u(t,p))$ is a $C^1$ function with respect to both arguments $(t,p)$. In broad descriptions, these conditions include certain smoothness of the
functions in Problem~(OCP1), satisfaction of the strict Legendre--Clebsch condition, uniqueness of the optimal control
minimizing the Hamiltonian, nonsingularity of the Jacobian of an
associated boundary-value problem, and boundedness of the symmetric solution of an associated Riccati ODE.

Fixing an increment $d  \in P$, the differentials 
$$
z_d(t,p_0) = \frac{\partial x}{\partial p}(t,p_0)d, \quad
v_d(t,p_0) = \frac{\partial u}{\partial p}(t,p_0)d ,
$$
satisfy a linear boundary value problem that contains only information obtained in the process of computing the unperturbed solution. The computations of these sensitivity differentials can also be performed by discretization methods
applied to the parametric optimal control problem; see B\"uskens \cite{Bueskens} and 
B\"uskens and Maurer \cite{Bueskens-Maurer}. The sensitivity differentials 
can be  conveniently used in the minimization of a master function defined on the Pareto front; see Section~\ref{sec:OPF}.

The above mentioned conditions for showing solution differentiability exclude optimal control problems
with control appearing linearly, since for this class of problems the strict Legendre-Clebsch condition does not hold.
  Here,  optimal controls are  combinations of bang-bang and singular arcs.
In case of finitely many switching times   and junction times with the boundary of a mixed control-state constraint
or a pure state constraint, one can set up a finite-dimensional  optimization problem, the Induced Optimization Problem, where the switching and junction times are optimized directly; see Maurer et al. \cite{MauBueKimKay2005} and
Osmolovskii and Maurer \cite{Osmolovskii-Maurer-book}.
 If second-order sufficient conditions hold for the Induced Optimization Problem (see \cite{Osmolovskii-Maurer-book}), 
one immediately obtains the result that the switching and junction times locally are differentiable functions
 of the parameter $p$. 

To our knowledge  extensions of these results on solution differentiability to optimal control problems with control 
and state delays can  not be found in the literature.

\section{An Algorithm For Optimization Over the Pareto Front}
\label{algorithm}

As discussed in Section~\ref{sec:soldiff}, the results~\cite[Theorem~3.1]{MauPes1994}~and~\cite[Theorem~5.1]{MauPes1995} lay the ground for devising and implementing numerical methods for
solving Problem~(OPF).  Bonnel and Kaya propose in~\cite{BonKay2010} a barrier method for convex bi-objective optimal control problems with pure control constraints. Their method relies on twice continuous differentiability of the solution (class $C^2$) in the weight $w$, using the result in~\cite[Theorem~3.1]{MauPes1994}.  

In this paper, we propose a bisection method also for the case of two objectives, which relies on the solution of Problem~(OCP$_w$) being of class $C^1$ w.r.t. the weight $w$, and thus taking the result in~\cite[Theorem~5.1]{MauPes1995} as a basis.  Although a mathematical justification of the applicability of our proposed method, i.e., solution differentiability, is given only for Problem~(OCP$_w$), the working of the method will also be illustrated on problems of more general class as in Problem~(OCP$_w$).

In the scalarized problem~(OCP$_w$) with two objectives ($r = 2$), by choosing $w_1 = w$, and $w_2=1 - w$, where $w\in[0,1]$, one can simply consider the single parameter $w$.

\subsection{Essential interval of weights}
\label{sec:ess_weights}

\begin{figure}
\begin{center}
\psfrag{pf}{Pareto front}
\psfrag{phi1}{$\varphi_1$}
\psfrag{phi2}{$\varphi_2$}
\psfrag{p1}{$(\varphi_1^*,\varphib_2)$}
\psfrag{p2}{$(\varphib_1,\varphi_2^*)$}
\psfrag{wf}{$w_f\,(\varphi_1 - \beta_1^*) = (1-w_f)\,(\varphi_2 - \beta_2^*)$}
\psfrag{w0}{$w_0\,(\varphi_1 - \beta_1^*) = (1-w_0)\,(\varphi_2 -
  \beta_2^*)$}
\psfrag{utop}{$(\beta_1^*,\beta_2^*)$}
\includegraphics[width=80mm]{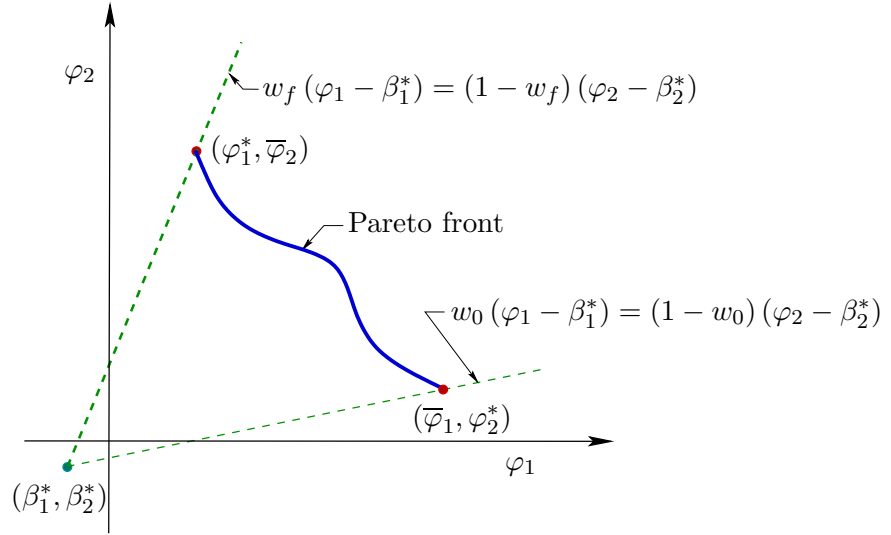}
\caption{\sf Determination of the essential subinterval of weights $[w_0,w_f]$~\cite{KayMau2014}.}
\end{center}
\label{fig:weights}
\end{figure}

With the Chebyshev scalarization, it would usually be enough for the weight $w$ to take values over a (smaller) subinterval $[w_0,w_f]\subset[0,1]$, with $w_0 > 0$ and $w_f<1$, for the generation of the whole front.  Figure~\ref{fig:weights} illustrates the geometry to compute the subinterval end-points, $w_0$ and $w_f$.  In the illustration, the points $(\varphi_1^*,\varphib_2)$ and $(\varphib_1,\varphi_2^*)$ represent the boundary of the Pareto front.  The equations of the ``rays'' which emanate from the utopia point $(\beta_1^*,\beta_2^*)$ and pass through the boundary points are also shown.  By substituting the boundary values of the Pareto curve into the respective equations, and solving each equation for $w_0$ and $w_f$ one simply gets
\begin{equation} \label{w0wf}
w_0 = \frac{(\varphi_2^* - \beta_2^*)}{(\varphib_1 - \beta_1^*) +
  (\varphi_2^* - \beta_2^*)}
\qquad\mbox{and}\qquad
w_f = \frac{(\varphib_2 - \beta_2^*)}{(\varphi_1^* - \beta_1^*) +
  (\varphib_2 - \beta_2^*)}\,.
\end{equation}

From the geometry depicted in Figure~\ref{fig:weights}, as also discussed in~\cite{KayMau2014}, one can deduce that with every $w\in[0,w_0]$ the solution of (OCP$_w$) will yield the same boundary point $(\varphib_1,\varphi_2^*)$ on the Pareto front.  Likewise with every $w\in[w_f,1]$ the same boundary point $(\varphi_1^*,\varphib_2)$ is generated.  This observation justifies the avoidance of the weights $w\in[0,w_0)\cup(w_f,1]$ in order not to keep getting the boundary points of the Pareto front, as otherwise one would end up wasting valuable computational effort and time.  

As a result of the above argument, the bisection method, implemented in the algorithm described in the next section, starts with the {\em essential interval} $[w_0,w_f]$ rather than $[0,1]$.  It is worth re-iterating that our main concern here, unlike in~\cite{KayMau2014}, is not really to construct the Pareto front, but rather do a search (in this case using the bisection method) over the Pareto front, at the same time avoiding the task of constructing the front, so as to find in some sense the best solution point in the Pareto front.

\subsection{Bisection method for solving Problem~(OPF)}

The problem of finding a best point in the Pareto front/set has now been transformed into a problem of finding best $w$, by virtue of the surjection from the set of weights to the set of Pareto minima furnished by Theorem~\ref{scalar_theorem}.  This has resulted in Problem~(OPF) and its concise form: Find some weight $w\in[w_0,w_f]$ such that the master objective function $\varphi_0(x^w, u^w, t_f^w)$ is minimized, where $(x^w, u^w, t_f^w)$ is found by solving (OCP$_w$) for that $w$.  For a simpler setting, it is helpful to define a function $F:[0,1]\to\dR_+$ representing the function we want to minimize over the Pareto front:
\begin{equation} \label{def:F}
F(w) := \varphi_0(x^w, u^w, t_f^w)\,, 
\end{equation}
such that $(x^w, u^w, t_f^w)$ solves (OCP$_w$).  In other words, an evaluation of the function $F(\cdot)$ at $w$ requires the solution of Problem~(OCP$_w$) with that $w$.  

Problem~(OPF) can now be re-written in an even more concise form as
\begin{equation} \label{OPF_concise}
\min_{w\in[w_0,w_f]}\ F(w)\,,
\end{equation}
where $F(\cdot)$ is evaluated as in~\eqref{def:F}.  In~\cite{BonKay2010}, a log-barrier method is proposed and implemented to solve~\eqref{OPF_concise}, with an underlying convex and smooth optimal control problem with no state constraints for which the solution can be assumed to be of class $C^2$, and so Newton-like methods are used with heuristic barrier parameter updates.  For the general form we have in Problem~(OCP), which is nonconvex and has state constraints, we assume that the solution is of class $C^1$.  As elaborated in Section~\ref{sec:soldiff}, under certain regularity conditions which can in many cases be checked, this assumption is guaranteed to hold.  Therefore we apply the bisection method~\cite{BurFai2011} as an effective and simple approach to solving~\eqref{OPF_concise} in the case of this paper.

Albeit elementary and standard, a statement of the optimality conditions in the fact below will be useful in formulating a computational algorithm later in this section. 
\begin{fact} \label{loc_min}
Consider the minimization problem in~\eqref{OPF_concise} with $F(\cdot)$ of class $C^1$.
\begin{itemize}
    \item[(a)] The interior point $w^*\in(w_0,w_f)$ is a {\em strict local minimizer} of $F(\cdot)$ if, and only if,
    \begin{equation} \label{Fermat}
    F'(w^*) = 0\,,
    \end{equation}
    and, for arbitrarily small $\varepsilon>0$,
    \begin{equation} \label{deriv_test}
    F'(w^* - \varepsilon) < 0\quad\mbox{and}\quad
    F'(w^* + \varepsilon) > 0\,.
    \end{equation}
    \item[(b)] The end point $w_0$ (resp.\ $w_f$) is a {\em strict local minimizer} of $F(\cdot)$ if, and only if, {\em either}
    \begin{itemize}
        \item[(i)] $F'(w_0) > 0$ (resp.\ $F'(w_f) < 0$) {\em or}
        \item[(ii)] $F'(w_0) = 0$ (resp.\ $F'(w_f) = 0$) and, for arbitrarily small $\varepsilon>0$, $F'(w_0 + \varepsilon) > 0$ (resp. $F'(w_f - \varepsilon) < 0$).
    \end{itemize} 
\end{itemize}
\end{fact}

\begin{remark}[Three Cases for the End Points of \boldmath{$[w_0,w_f]$}]  \label{rem:3cases} \rm
We will apply the bisection method starting with the essential interval $[w_0,w_f]$.  Before introducing the pertaining algorithm, we consider below the cases for the end points of this interval.
\begin{description}
    \item[Case I.] $F'(w_0)F'(w_f) < 0$\,:  Since $F'(\cdot)$ is assumed to be continuous the bisection method is guaranteed to find a numerical solution to~\eqref{Fermat} by the intermediate value theorem.  Condition~\eqref{deriv_test} needs to be check to see if $w^*$ is a strict local minimizer.
    \item[Case II.] $F'(w_0)F'(w_f) > 0$\,:  By the conditions in~Fact~\ref{loc_min}(b)(i), at least one of $w_0$ and $w_f$ is a strict local minimizer.
    \item[Case III.] $F'(w_0)F'(w_f) = 0$\,:  If one of the inequalities in~Fact~\ref{loc_min}(b)(ii) is satisfied, then $w_0$ or $w_f$ is a strict local minimizer.  It is possible that both $w_0$ and $w_f$ are, or only one or neither is, a local minimizer.
\end{description}
In Case~I, the bisection method starts with the interval $[w_0,w_f]$ and terminates with an approximate solution in the interior of the interval.  In Case~II, a local minimum is found immediately, and so in principle there is no need to do a further search.  In Case~III, however, the conclusion might be that neither $w_0$ nor $w_f$ is a strict local minimizer, in which case it would be necessary to start the bisection method with a subinterval of $[w_0,w_f]$, and consider Cases~I--III again.
\proofbox
\end{remark}

\begin{remark}  \rm
In any of the scenarios elaborated in Remark~\ref{rem:3cases}, consideration of another subinterval of $[w_0,w_f]$ might as well yield a better (lower-value) solution, since the problem is nonconvex and we can only hope to get a locally optimal solution.  In our approach here, however, we do not endeavour to obtain a global minimum.  As a result of our discussion in Remark~\ref{rem:3cases}, we will consider only Case~I, which clearly prompts us to use the bisection method directly.  As suggested above, in the event of Case~III not yielding a solution, the new subinterval could be chosen in such a way that one would fall into Case~I.
\proofbox
\end{remark}

The derivative of $F(\cdot)$ is defined at the end points of the interval $[w_0,w_f]$ as one-sided limits,
\[
F'(w_0) := \lim_{\delta\to 0^+} \frac{F(w_0+\delta) - F(w_0)}{\delta} \quad\mbox{and}\quad
F'(w_f) := \lim_{\delta\to 0^-} \frac{F(w_f+\delta) - F(w_f)}{\delta}\,, \]
and in the interior, i.e., for $w\in (w_0,w_f)$, as
\[
F'(w) := \lim_{\delta\to 0} \frac{F(w+\delta) - F(w)}{\delta}\,, 
\]
where $F(\cdot)$ is evaluated as in~\eqref{def:F}.  In computations, we will use the forward, and backward, finite difference approximations of $F'(\cdot)$.  Namely, for some small $\delta > 0$, we will set
\begin{equation} \label{difference}
F'(w) \approx 
\left\{\begin{array}{rl}
\dfrac{F(w+\delta) - F(w)}{\delta}\,, &\ \mbox{ if } w\in[w_0,w_f-\delta)\,, \\[3mm] 
\dfrac{F(w) - F(w-\delta)}{\delta}\,, &\ \mbox{ if } w\in[w_f-\delta,w_f]\,. 
\end{array} \right.
\end{equation}
The {\em step} $\delta$ in the difference approximation formula \eqref{difference} is small for an accurate estimation of the derivative but not too small in order not to divide one very small number by another and cause numerical instabilities.

In what follows we provide an algorithm to solve Problem~(OPF).  The algorithm first finds the essential interval $[w_0,w_f]$, computes the signs of $F'(w_0)$ and $F'(w_f)$ and checks the cases I--III in Remark~\ref{rem:3cases}, and then if $F'(w_0)F'(w_f) < 0$ it uses the bisection method, to find a numerical solution to Problem~(OPF).

\begin{algorithm} \label{algo1} \rm \ 
\begin{description}
\item[Step {\boldmath$0.0$}] ({\em Initialization}) \ \ Choose utopia parameters,\ \ $\eta_1,\eta_2>0$, a small numerical differentiation step $\delta > 0$, a stopping tolerance $\epsilon > 0$, and a maximum number of iterations $k_{\max}$\,.\ \  Set\ \ $k:=1$.
\item[Step {\boldmath$0.1$}] ({\em Boundary points of the front})\ \ Solve \eqref{single_obj} to get\ \ $(\xb^i,\ub^i,\tb_f^i)$, $i=1,2$.\ \ Set \\[2mm]
$(\varphi_1^*,\varphib_2) := (\varphi_1(\xb^1(\tb_f^1),\tb_f^1), \varphi_2(\xb^1(\tb_f^1),\tb_f^1))$,\ \ $(\varphib_1,\varphi_2^*) := (\varphi_1(\xb^2(\tb_f^2),\tb_f^2), \varphi_2(\xb^2(\tb_f^2),\tb_f^2))$\,.
\item[Step {\boldmath$0.2$}] ({\em Utopia point})\ \ Set\ \ $\beta^*
  := (\beta_1^*,\beta_2^*)$ with\ \ $\beta^*_i := \varphi_i^* -
  \eta_i$,\ \ $i=1,2$.
\item[Step {\boldmath$0.3$}] ({\em Essential interval})\ \ Determine the subinterval\ \ $[w_0,w_f]\subset[0,1]$\ \ using \eqref{w0wf}.
\item[Step {\boldmath$0.4$}] ({\em Signs at end points})\ \ Compute $F'(w_0)$ and $F'(w_f)$ using \eqref{difference}, with $F(\cdot)$ evaluated as in~\eqref{def:F}.
    \begin{itemize}
        \item If Fact~\ref{loc_min}(b)(i) or (ii) is satisfied then $w^* = w_0$ or $w^* = w_f$ appropriately; STOP.
        \item If $F'(w_0)F'(w_f) = 0$ and neither of the inequalities in~Fact~\ref{loc_min}(b)(ii) is satisfied then declare ``Algorithm failed. Change the interval $[w_0,w_f]$.'' and STOP.
    \end{itemize}
Let $a := w_0$ and $b := w_f$.
\item[Step {\boldmath$k.1$}] ({\em Bisection})\ \ Find the midpoint $c := a + (b-a)/2$ of the interval $[a,b]$.
\item[Step {\boldmath$k.2$}] ({\em Stopping criterion})\ \ Compute $F'(c)$  using \eqref{difference}, with $F(\cdot)$ evaluated as in~\eqref{def:F}.
    \begin{itemize}
        \item If $F'(c) = 0$ or $(b-a)/2 < \epsilon$ then set $w^* = c$ and STOP.
        \item If $k = k_{\max}$ then declare ``Maximum number of iterations exceeded.'' and STOP.
    \end{itemize}
\item[Step {\boldmath$k.3$}] ({\em New subinterval}) \ \ Set\ \ $k:=k+1$\,.  If $F'(a)F'(c) > 0$ then update the subinterval as $[a,b] := [c,b]$; otherwise, set $[a,b] := [a,c]$.
GO TO Step~$k.1$.
\end{description}
\end{algorithm}

\section{Numerical Examples}
\label{examples}

In this section, we illustrate the working of Algorithm~\ref{algo1} on two optimal control problems, one involving an electric circuit in Section~\ref{sec:Rayleigh} and the other a tuberculosis~(TB) epidemic in Section~\ref{sec:TB}.

In computations, we use direct discretization of optimal control problems for which convergence theory has been an active topic of research in the literature (see for example \cite{AltBaiLemGer2013, Betts2020, DonHag2001, DonHagMal2000, DonHagVel2000, PieScaVel2018}, and see \cite{KayMau2014} for additional references and discussion).

We employ the scalarize--discretize--then--optimize approach that was previously used in~\cite{KayMau2014}.  Under this approach, one first scalarizes the multi-objective problem in the infinite-dimensional space, and then discretizes the scalarized problem directly and applies a usually large-scale finite-dimensional optimization method to find a discrete approximate solution of the scalarized problem. By the existing theory of discretization mentioned above, under certain assumptions, the discrete approximate solution converges to a solution of the continuous-time scalarization of the original problem, yielding a Pareto minimum of the original problem.  When possible, we will also check {\em a posteriori} to see if the necessary optimality conditions are satisfied by an accurate-enough numerical solution. 

In Step~$0.4$ of Algorithm~\ref{algo1}, a direct discretization of Problem~(OCP$_w$), for example employing a Runge--Kutta scheme, such as Euler's method or the Trapezoidal rule, is solved by using Ipopt, version 3.12.13, four times.  In Step~$k.2$, Problem~(OCP$_w$) is solved in a similar way two times.   Ipopt is a popular optimization software based on an interior point method; see~\cite{WacBie2006}. We use AMPL \cite{AMPL} as an optimization modelling language, which employs Ipopt as a solver.

\subsection{Example: Tunnel-diode oscillator (Rayleigh problem)}
\label{sec:Rayleigh}

The tunnel-diode oscillator problem, also referred to as the Rayleigh problem in the literature, involves dynamics represented by the following differential equations.
\begin{eqnarray*}
&& \dot{x}_1(t) = x_2(t)\,,  \\[2mm]
&& \dot{x}_2(t) = -x_1(t) + x_2(t)\,(1.4 - 0.14\,x_2^2(t)) + 4\,u(t)\,, \mbox{\ \ for a.e. } t\in[0,t_f]\,,
\end{eqnarray*}
where the state variable $x_1(t)$ denotes electric current, and the control variable $u(t)$ stands for a suitable transformation of the voltage at a generator, both at time $t\in[0,t_f]$---see~\cite{MauObe2002} for a detailed exposition of the problem.  In this particular instance of the problem, the initial and terminal values of the state variables are specified as
\[
(x_1(0), x_2(0)) = (-5, -5)\quad\mbox{and}\quad (x_1(t_f), x_2(t_f)) = (0,0)\,,
\]
and the dynamics are subject to constraints on the control variable such that
\[
-1\le u(t)\le 1\,, \mbox{\ \ for a.e. } t\in[0,t_f]\,. 
\]
The optimal control problem is posed as a bi-objective problem with
\[
\min\ \left[ t_f\,,\ \ \int_0^{t_f} \left(x_1^2(t) + u^2(t)\right) dt\right]\,,
\]
where the competing objectives are the minimization of the final time $t_f$ and the minimization of the sum of the square $L^2$-norms, or in some sense the magnitudes, of the current and the generator voltage. Define a new state variable $x_3$ such that
\[
\dot{x}_3(t) = x_1^2(t) + u^2(t)\,, \mbox{\ \ for a.e. } t\in[0,t_f]\,,\ \ x_3(0) = 0\,.
\]
Then the two objective functionals as in Problem~(OCP), or Problem~(OCP$_w$), can be expressed as
\[
\varphi_1(x(t_f),t_f) = t_f\quad\mbox{and}\quad 
\varphi_2(x(t_f),t_f) = x_3(t_f)\,.
\]
As we have stated above, the bi-objective Rayleigh problem is in the same form as Problem~(OCP) and, in particular, Problem~(OCP$_w$).  The decision maker's objective for this problem will be to minimize a weighted distance to the origin of the value space.  We choose
\[
\varphi_0(x^w,u^w,t_f^w) := 100\,\varphi_1^2(x^w(t_f),t_f^w) + \varphi_2^2(x^w(t_f),t_f^w)\,,
\]
where the scaling multiplier 100 is used to make the orders of magnitudes of $\varphi_1$ and $\varphi_2$ the same.  We aim to solve Problem~(OPF), to determine a scalar $w\in(0,1)$ with $w_1 := w$ and $w_2 := 1-w$ that results in the best Pareto solution in the sense that $\varphi_0(\cdot,\cdot,\cdot)$ is minimized, subject to the solution of Problem~(OCP$_w$).

In~\cite{MauObe2002}, Maurer and Oberle numerically illustrate that an optimal solution does not exist for the single objective problem minimizing the quadratic functional $\varphi_2(x(t_f),t_f)$, in that $t_f$ tends to infinity.  They carry out a numerical test for checking the second-order sufficient conditions (SSC) of optimality and show that the test fails to confirm the SSC.  Therefore, we will impose a bound on the terminal time, namely set $t_f \le 5$.  On the other hand, they illustrate also in~\cite{MauObe2002} that for certain instances of the weighted-sum problem, the SSC of optimality are satisfied.

Problem~(OCP$_w$) can now explicitly be written for the Rayleigh problem as
\[
\left\{\begin{array}{rll}
\ds\min_{{\alpha\ge0}\atop{x(\cdot),u(\cdot),t_f}} & \ \ \alpha & \\[4mm]
\mbox{subject to} 
    & \ \ \dot{x}_1(t) = x_2(t)\,, & x_1(0) = -5\,,\ x_1(t_f) = 0\,, \\[2mm]
    & \ \ \dot{x}_2(t) = -x_1(t) + x_2(t)\,(1.4 - 0.14\,x_2^2(t)) + 4\,u(t)\,, & x_2(0) = -5\,,\ x_2(t_f) = 0\,, \\[2mm]
    & \ \ \dot{x}_3(t) = x_1^2(t) + u^2(t)\,, & x_3(0) = 0\,, \\[2mm]
    & \ \ -1\le u(t)\le 1\,, \mbox{\ \ for a.e. } t\in[0,t_f]\,,\ \ t_f \le 5\,, & \\[2mm]
    & \ \ w\,(t_f - \beta_1^*) \le\alpha\,, & \\[2mm]
	& \ \ (1-w)\,(x_3(t_f) - \beta_2^*) \le\alpha\,. & 
\end{array} \right.
\]

The Hamiltonian $H:\dR^3\times\dR\times\dR^3\to\dR$ for this problem simply is
\[
H(x,u,\lambda) := 
\lambda_1 x_2 + \lambda_2\left[(-x_1 + x_2\,(1.4 - 0.14\,x_2^2) + 4 u\right] + 
\lambda_3 (x_1^2 + u^2)\,,
\]
where $\lambda(t) := (\lambda_1(t),\lambda_2(t),\lambda_3(t))\in\dR^3$ is referred to as the adjoint variable vector.  Using the convenient notation $H[t] := H(x(t),u(t),\lambda(t))$, suppose that
\begin{subequations}
\begin{eqnarray}
&& \dot{\lambda}_1(t) := -H_{x_1}[t] = \lambda_2(t) - 2\lambda_3(t) x_1(t)\,, \label{adj1} \\[1mm]
&& \dot{\lambda}_2(t) := -H_{x_2}[t] = -\lambda_1(t) - \lambda_2(t) (1.4 - 0.42 x_2^2(t))\,, \\[1mm]
&& \dot{\lambda}_3(t) := -H_{x_3}[t] = 0\,,  \label{adj3}
\end{eqnarray}
\end{subequations}
for all $t\in[0,t_f]$, with certain transversality conditions as required by the maximum principle.  In \eqref{adj1}--\eqref{adj3}, $H_{x_i} := \partial H/\partial x_i$, $i = 1,2,3$.  We will not go into the details of these (boundary) conditions here.  However we note that $\lambda_3(t) = \overline{\lambda}_3$, a constant, for all $t\in[0,t_f]$.  Then the maximum principle states that if $(x,u,t_f)$ is an optimal solution triplet then there exists a continuous function $\lambda(\cdot)$ satisfying \eqref{adj1}--\eqref{adj3}, along with certain transversality conditions, such that $\lambda(t) \neq 0$, for all $t\in[0,t_f]$, and
\begin{equation}  \label{gen_opt_cont}
u(t) = \argmin_{v\in[-1,1]} H(x(t),v,\lambda(t)) = 
    \argmin_{v\in[-1,1]} \left(4\lambda_2(t) v + \lambda_3(t) v^2\right).
\end{equation}
for a.e. $t\in[0,t_f]$.  If $w = 1$, then the problem is a single-objective one, referred to as a time-optimal control problem, and the condition \eqref{gen_opt_cont} reduces to
\begin{equation*}  \label{toc_opt_cont}
u(t) = \argmin_{v\in[-1,1]} \lambda_2(t) v\,,
\end{equation*}
resulting in
\begin{equation}  \label{toc_opt_cont1}
u^w(t) = \left\{\begin{array}{rl}
	1\,, &\ \ \mbox{if\ \ } \lambda^w_2(t) < 0\,, \\[1mm]
	-1\,, &\ \ \mbox{if\ \ } \lambda^w_2(t) > 0\,, \\[1mm]
	\mbox{undetermined}\,, &\ \ \mbox{if\ \ } \lambda^w_2(t) = 0\,,
\end{array} \right.
\end{equation}
for a.e. $t\in[0,t_f]$.  By the discussion given in Section~\ref{sec:ess_weights} (also see~\cite{KayMau2014}), $u^w(t)$ given in \eqref{toc_opt_cont1} is the same for all $w\in[w_f, 1]$.  Recall that if one does not have $\lambda^w_2(t) = 0$ for all $[t',t'']\subset[0,t_f]$, where $t' < t''$, then $u^w(t)$ in \eqref{toc_opt_cont1} is referred to as optimal control of {\em bang--bang} type.  We assume (and therefore will numerically double-check) that the optimal control for the particular instance of the problem is of bang--bang type.

\begin{figure}[t]
\begin{minipage}{80mm}
\begin{center}
\includegraphics[width=80mm]{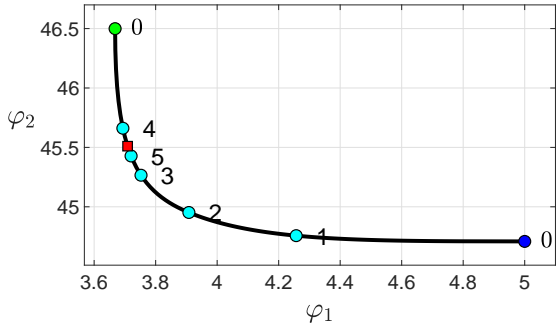} \\[3mm]
{\sf\small (a) Pareto front, and iterations of Algorithm~\ref{algo1}: Master solution is depicted by a (red) square and iterates by (light blue) circles.}
\end{center}
\end{minipage}
\begin{minipage}{80mm}
\begin{center}
\includegraphics[width=80mm]{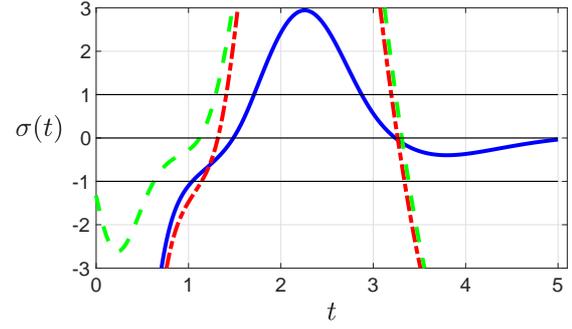} \\[3mm]
{\sf\small (b) Switching function as defined in~\eqref{switch_fun}. \\[6mm] \ }
\end{center}
\end{minipage}
\\[10mm]
\begin{minipage}{80mm}
\begin{center}
\psfrag{singular switching curve}{singular control switching curve}
\includegraphics[width=90mm]{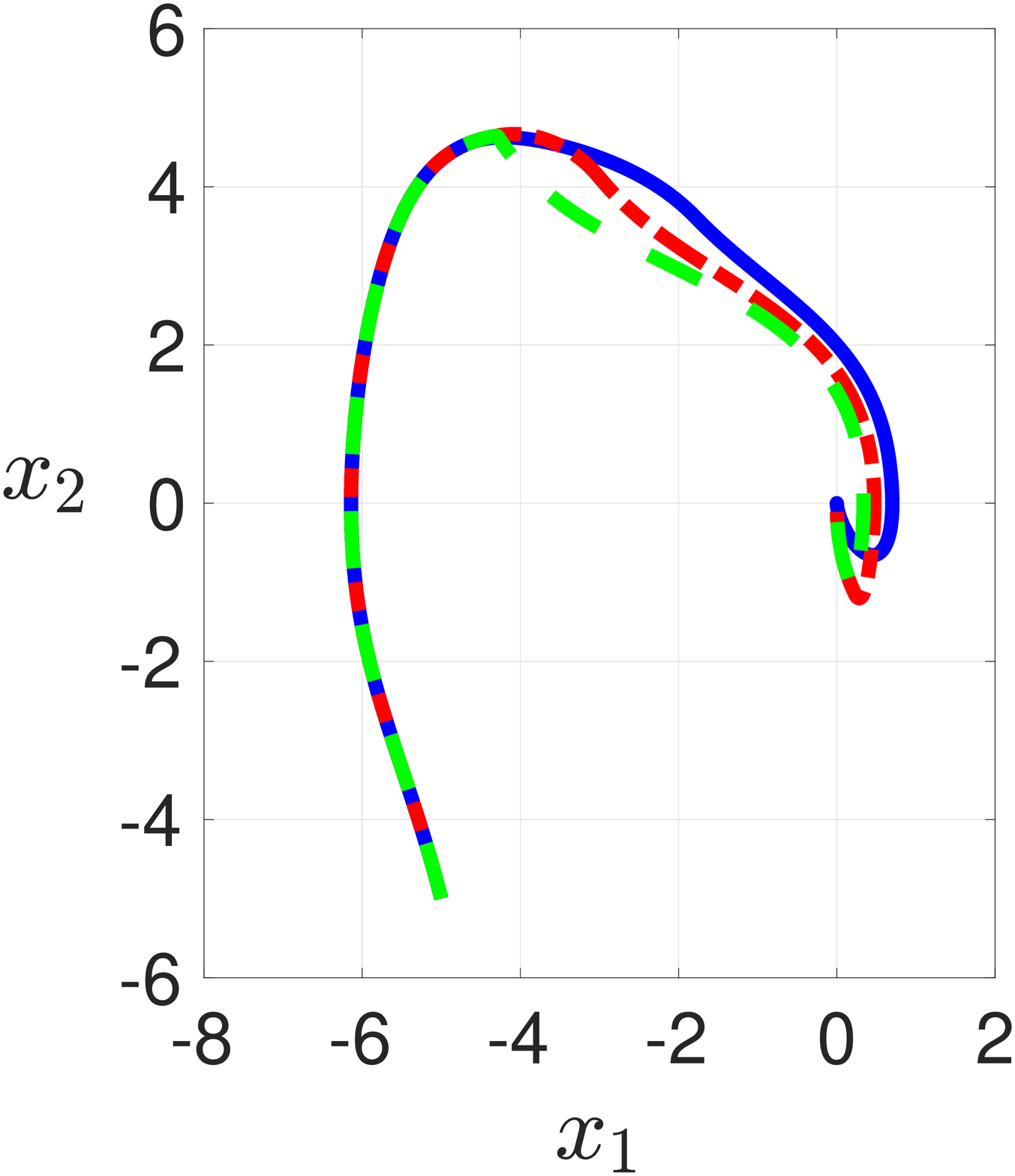} \\[3mm]
{\sf\small (c) Phase plane trajectories. \\[1mm] \ }
\end{center}
\end{minipage}
\begin{minipage}{80mm}
\begin{center}
\includegraphics[width=80mm]{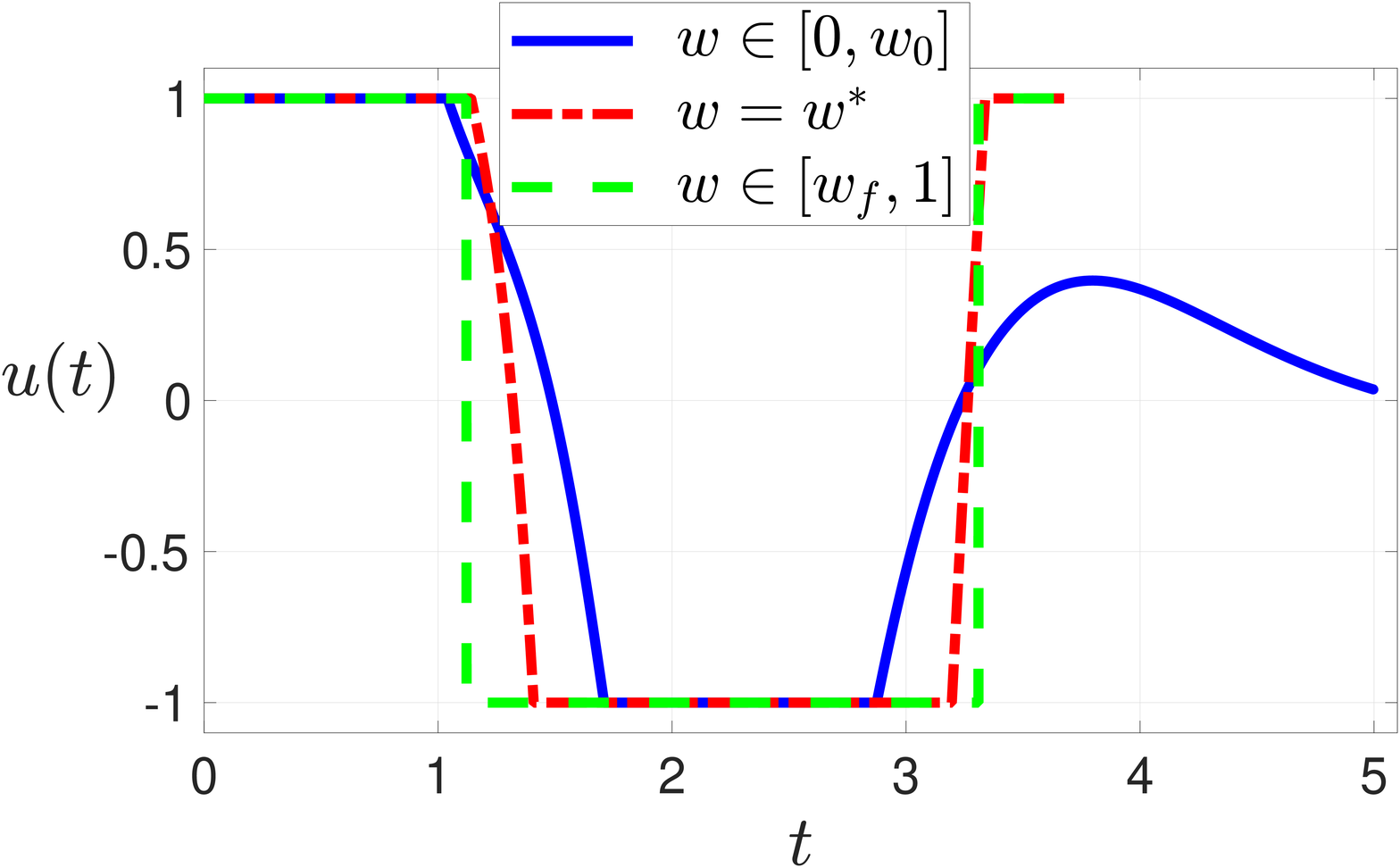} \\[3mm]
{\sf\small (d) Control variable.} 
\end{center}
\end{minipage}
\
\caption{\sf Rayleigh problem---Boundary Pareto solutions, corresponding to $w_0 = 0.8994$ and $w_f = 0.9269$, are shown with (blue) solid curves and (green) dashed curves, respectively.  Master Pareto solution, corresponding to $w^* = 0.9247$, is shown with dashed-and-dotted (red) curves.}
\label{fig:Rayleigh}
\end{figure}

The optimality condition~\eqref{gen_opt_cont} can be shown to yield, for any given $w\in[w_0,w_f)$,
\begin{equation}  \label{gen_opt_cont1}
u^w(t) = \left\{\begin{array}{rl}
	1\,, &\ \ \mbox{if\ \ } 2\lambda^w_2(t) < -\overline{\lambda}^w_3\,, \\[1mm]
	-2\lambda^w_2(t)/\overline{\lambda}^w_3\,, &\ \ \mbox{if\ \ } -\overline{\lambda}^w_3 \le 2\lambda^w_2(t) \le \overline{\lambda}^w_3\,, \\[1mm]
	-1\,, &\ \ \mbox{if\ \ } 2\lambda^w_2(t) > \overline{\lambda}^w_3\,, \\[1mm]
\end{array} \right.
\end{equation}
for all $t\in[0,t_f]$, provided $\overline{\lambda}^w_3 \neq 0$.  Again by virtue of the discussion in Section~\ref{sec:ess_weights}, $u^w(t)$ in \eqref{gen_opt_cont1} is the same for all $w\in[0,w_0]$.  We define the {\em switching function} as
\begin{equation}  \label{switch_fun}
\sigma^w(t) := \left\{\begin{array}{rl}
2\,\lambda^w_2(t)/\overline{\lambda}^w_3\,, &\ \ \mbox{if\ \ } 0\le w < w_f\,, \\[1mm]
16\,\lambda^w_2(t)\,, &\ \ \mbox{if\ \ } w_f\le w \le 1\,.
\end{array} \right.
\end{equation}
The constant coefficients 2 and 16 above are used for scaling purposes, so that the graphs in Figure~\ref{fig:Rayleigh}(b) can be viewed more easily. Now, using \eqref{switch_fun}, we can summarize and combine the expressions for the optimal control in \eqref{toc_opt_cont1} and \eqref{gen_opt_cont1} as follows.
\begin{equation}  \label{opt_cont}
u^w(t) = \left\{\begin{array}{rl}
\left\{\begin{array}{rl}
	1\,, &\ \ \mbox{if\ \ } \sigma^w(t) < -1\, \\[1mm]
	-2\lambda^w_2(t)/\overline{\lambda}^w_3\,, &\ \ \mbox{if\ \ } -1 \le \sigma^w(t) \le 1\,, \\[1mm]
	-1\,, &\ \ \mbox{if\ \ } \sigma^w(t) > 1\,. \\[1mm]
\end{array} \right\}, &\ \ \mbox{if\ \ } 0\le w < w_f\,, \\[7mm]
\left\{\begin{array}{rl}
	1\,, &\ \ \mbox{if\ \ } \sigma^w(t) < 0\,, \\[1mm]
	-1\,, &\ \ \mbox{if\ \ } \sigma^w(t) > 0\,. \\[1mm]
\end{array} \right\}, &\ \ \mbox{if\ \ } w_f\le w \le 1\,.
\end{array} \right.
\end{equation}
As to why $\sigma^w(\cdot)$ is referred to as the switching function should now be more clear from~\eqref{opt_cont}: the value of $\sigma^w(\cdot)$ determines when to switch from one case of the control function $u^w(\cdot)$ to another.

For Problem~(OCP$_w$) written for the Rayleigh problem above, we have chosen the utopia vector as $(\beta_1^*, \beta_2^*) = (0, 0)$, since $\varphi_i(x(t_f),t_f) > 0$, for $i=1,2$.  Figure~\ref{fig:Rayleigh}(a) depicts the Pareto front for the instance of the multi-objective Rayleigh problem we consider here.  It also displays the iterations of Algorithm~\ref{algo1}.  The Rayleigh problem is discretized using the trapezoidal rule, the number of grid points is set to be $N = 5000$, and the Ipopt's tolerance to $10^{-10}$, so as to get solutions for $w$ accurate at least up to four decimal places~(dp).

The essential interval is found to be $[w_0,w_f] = [0.8994, 0.9269]$, with
\[
(\varphi_1^{w_0},\varphi_2^{w_0}) = (5.000, 44.71)\quad\mbox{and}\quad
(\varphi_1^{w_f},\varphi_2^{w_f}) = (3.668, 46.50)\,,
\]
correct to four significant figures, where $\varphi_i^{w} := \varphi_i(x^{w}(t_f),t_f^{w})$, $i = 1, 2$, with $w = w_0$ or $w_f$, or as will be the case below, $w = w^*$.
Optimization over the Pareto front results in $w^* = 0.9247$, after 14 iterations of Algorithm~\ref{algo1}, yielding 
\[
\varphi_0^{w^*} = 58.71 \quad\mbox{and}\quad
(\varphi_1^{w^*},\varphi_2^{w^*}) = (3.709,45.51)\,.
\]
If there is a need to save the computational resources further, the algorithm can be asked to yield a less accurate result, say correct to three dp, which then yields $w^* = 0.925$ in eight iterations with $(\varphi_1^{w^*},\varphi_2^{w^*}) = (3.71,45.5)$.  In Figure~\ref{fig:Rayleigh}(a) only five iterations are displayed (labels 1--5 appearing to the right of each iteration) for clarity in viewing.  The Pareto (master) solution with $w = w^*$ is represented by a square.

The numerical Pareto-optimal state and control variable solutions are presented in~Figures~\ref{fig:Rayleigh}(c)--(d) for $w = w_0, w^*, w_f$.  One of the boundary Pareto-optimal solutions is shown using solid (blue) curves for $w = w_0$, which is the same solution for all $w\in[0,w_0]$, as previously discussed in~Section~\ref{sec:ess_weights}.  On the other hand, the other boundary Pareto-optimal solution for $w = w_f$, which holds for all $w\in[w_f,1]$, is shown using dashed (green) curves.  The latter is nothing but a time-optimal control solution for the Rayleigh problem (a solution with the smallest $t_f$), resulting in a bang--bang type function with the sequence of values $\{1,-1,1\}$, namely with two switchings.  The master Pareto solution is given for $w = w^*$ using dashed-and-dotted (red) curves.

The switching function $\sigma^w(\cdot)$ plotted in~Figure~\ref{fig:Rayleigh}(b) by using~\eqref{switch_fun} (recall that discrete approximations of $\lambda^w_2(t)$ and $\lambda^w_3(t)$ can readily be obtained from AMPL) furnishes the means to verify the optimality condition for $u^w(\cdot)$ expressed in~\eqref{opt_cont}.  It is evident from the dashed (green) plot of the switching function that, for $w\in[w_f,1]$, when $\sigma^w(\cdot)$ crosses the time axis there is a jump (from $1$ to $-1$ or vice versa) in the value of the corresponding $u^w(\cdot)$ plot.  Likewise, for $w\in[0,w_0]$ and for $w = w^*\in[w_0,w_f)$, whenever $\sigma^w(\cdot)$ crosses one of the lines $\sigma^w(t) = 1$ and $\sigma^w(t) = -1$ (shown by two black lines in~Figure~\ref{fig:Rayleigh}(b) for convenience) the expression for the control function $u^w(\cdot)$ switches from one case in~\eqref{opt_cont} to another, as required.

\subsection{Example: Compartmental model for tuberculosis}
\label{sec:TB}

In 2020 and 2021, tuberculosis (TB) was the second leading cause of death from an infectious disease worldwide after COVID-19~\cite{TB_WHO2021}.  Active TB refers to disease that occurs in someone infected with Mycobacterium tuberculosis.  It is characterized by signs or symptoms of active disease, or both, and is distinct from latent tuberculosis infection, which occurs without signs or symptoms of active disease. Only individuals with active TB
can transmit the infection. Many people with active TB do not experience typical TB symptoms in the early stages of the disease. These individuals are unlikely to seek care early, and may not be properly diagnosed when seeking care. Delays to diagnosis of active TB present a major obstacle to the control of a TB epidemic, it may worsen the disease, increase the risk of death and enhance
tuberculosis transmission to the community. Both patient and the health system may be responsible for the diagnosis delay.

We study the control model with control and state delays presented in Silva et al. \cite{Silva-Maurer-Torres}.  In this model, reinfection and post-exposure interventions for tuberculosis are considered. The population is divided into five categories (compartments) (\textrm{i.e.}, the control system  has five state variables):

%%============================================================

\begin{tabular}{lcl}
$S$ &:& susceptible individuals, \\
$L_1$ &:& early latent individuals, recently infected (less than two years), \\
$I$ &:& infectious individuals, who have active TB, \\
$L_2$ &:& persistent latent individuals, \\
$R$ &:& recovered individuals, \\
$N$ &:& total population $N = S + L_1 + I + L_2 + R$ , assumed constant.
\end{tabular}

\vspace{1mm}
The model has two control variables and three delays:

\begin{tabular}{lcl}
$u_1$ &:& effort on early detection and treatment of recently infected individuals $L_1$, \\
$d_{u_1}$ &:& delay on the diagnosis of latent TB,  and commencement of latent TB treatment, \\
$u_2$ &:& chemotherapy or post-exposure vaccine to persistent latent individuals $L_2$, \\
$d_{u_2}$ &:& delay in the prophylactic treatment of persistent latent $L_2$, \\
$d_I$ &:& delay in $I$, i.e., delay in diagnosis.
\end{tabular}

\vspace{1mm}
The dynamical system is given by
% \begin{eqnarray}
% \label{eq:S}
% \dot{S}(t) &=& \mu N - \beta I(t) S(t) / N - \mu S(t)
% \end{eqnarray}
%%=====================================================
\begin{equation}
\label{dynamics}
\left\{
\begin{aligned}
\dot{S}(t) &= \mu N - \frac{\beta}{N} I(t) S(t) - \mu S(t),  \\[1mm]
\dot{L}_1(t) &= \frac{\beta}{N} I(t)\left( S(t)
+ \sigma L_2(t) + \sigma_R R(t)\right) - \left(\delta + \tau_1 + \epsilon_1 u_1(t-d_{u_1}) + \mu\right)L_1(t), \\[1mm]
%%-----------------------------------------------------
\dot{I}(t) &= \phi \,\delta \, L_1(t) + \omega L_2(t) + \omega_R R(t)
- \tau_0 I(t-d_I)\ + \mu I(t), \\[1mm]
%%----------------------------------------------
\dot{L}_2 (t) &= (1 - \phi) \delta L_1(t) - \sigma \frac{\beta}{N} I(t) L_2(t)
- (\omega + \epsilon_2 u_2(t-d_{u_2}) + \tau_2 + \mu)L_2(t).
%%---------------------------------
% \\
%%---------------------------------------------------------
% \dot{R}(t) &= \tau_0 I(t) +  (\tau_1 + \epsilon_1 u_1(t)) L_1(t)
% + (\tau_2 +\epsilon_2 u_2(t)) L_2(t) - \sigma_R \frac{\beta}{N} I(t) R(t)\\
% &\quad \quad \quad - (\omega_R + \mu)R(t) \, .
\end{aligned}
\right.
\end{equation}
%%======================================================
The recovered population is defined by % $R := N - S-L_1 - I - L_2$ 
\begin{equation}  \label{R}
R(t) := N - S(t) - L_1(t) - I(t) - L_2(t)\,,
\end{equation}
with $N=30000$.  
The system and delay parameters in the model~\eqref{dynamics} along with their values are listed in Table~\ref{table:parameters}. 
% {The delays are set as $d_{u_1} = d_{u_2} = 0.2$ and $d_I = 0.1$.}
In view of the delays the initial conditions and functions are:
\begin{equation}
\label{initial-condition}
\hspace{-0mm}
\begin{array}{l}
S(0) = 76\,N/120, \ L_1(0) = 36\,N/120, \ L_2(0) = 2\,N/120 , \ R(0) = N /120, \\[2mm]
I(t) = 5\,N/120 \quad \mbox{for} \ -\!d_I \leq t \leq 0 , \quad
 u_k(t) = 0 \quad \mbox{for} \ -\!d_{u_k} \leq t < 0 , \quad (k=1,2).
\end{array}
\end{equation}
%%============================================================
%% ------------------------------------------------
\begin{table}
\centering
\begin{tabular}{|l|l|l|}
\hline
{\small{Symbol}} & {\small{Description}}  & {\small{Value}} \\
\hline
{\small{$\beta$}} & {\small{Transmission coefficient}}  & {\small{variable}}\\
{\small{$\mu$}} & {\small{Death and birth rate}}  & {\small{$1/70 \, yr^{-1}$}}\\
{\small{$\delta$}} & {\small{Rate at which individuals leave $L_1$}}  & {\small{$12 \, yr^{-1}$}}\\
{\small{$\phi$}} & {\small{Proportion of individuals going to $I$}}  & {\small{$0.05$}}\\
{\small{$\omega$}} & {\small{Endogenous reactivation rate for persistent latent infections}}  & {\small{$0.0002 \, yr^{-1}$}}\\
{\small{$\omega_R$}} & {\small{Endogenous reactivation rate for treated individuals}}   &{\small{$0.00002 \, yr^{-1}$}}\\
{\small{$\sigma$}} & {\small{Factor reducing the risk of infection as a result of acquired}}  & \\
& {\small{immunity to a previous infection for $L_2$}} & {\small{$0.25$}} \\
{\small{$\sigma_R$}} & {\small{Rate of exogenous reinfection of treated patients}}  & {\small{0.25}} \\
{\small{$\tau_0$}} & {\small{Rate of recovery under treatment of active TB}}  &  {\small{$2 \, yr^{-1}$}}\\
{\small{$\tau_1$}} & {\small{Rate of recovery under treatment of early latent individuals $L_1$}}  &  {\small{$2 \, yr^{-1}$}}\\
{\small{$\tau_2$}} & {\small{Rate of recovery under treatment of persistent latent individuals $L_2$}}  &  {\small{$1 \, yr^{-1}$}}\\
{\small{$N$}} & {\small{Total population}} & {\small{$30,000$}} \\
{\small{$\epsilon_1$}} & {\small{Efficacy of treatment of early latent $L_1$}} & {\small{$0.5$}} \\
{\small{$\epsilon_2$}} & {\small{Efficacy of treatment of persistent latent TB $L_2$}} & {\small{$0.5$}} \\
{\small{$t_f$}} & {\small{Total simulation duration}} & {\small{$5$ $years$}} \\
{\small{$d_I$}} & {\small{delay in the diagnosis of $I$}} & {\small{$0.1$ $years$}} \\
{\small{$d_{u_1}$}} & {\small{delay in the diagnosis of early latent individuals $L_1$}} & {\small{$0.2$ $years$}} \\
{\small{$d_{u_2}$}} & {\small{delay in the prophylactic treatment of persistent latent individuals $L_2$}} & {\small{$0.2$ $years$}} \\
%{\small{$W_1$}} & {\small{Weight constant of control $u_1$}} & {\small{$500$}}\\
%{\small{$W_2$}} & {\small{Weight constant of control $u_2$}} & {\small{$50$}}\\
\hline
\end{tabular}
\caption{Parameter values for the TB control model.}
\label{table:parameters}
\end{table}
%%=================================================================
The control constraints are given by
\begin{equation}
\label{control-constraint}
0 \leq u_k(t) \leq 1\,, \quad \forall  t\in[0, t_f] \,,\quad (k=1,2).
\end{equation}
We consider the following parametric objective functional with control weights $a_1, a_2 \geq 0$:
\begin{equation}
\label{objective-general}
\int\limits_0^{t_f} (I(t) + L_2(t) + a_1 u_1(t) + a_2 u_2(t)) \, dt \,.
\end{equation}
Depending on the priorities, the weights $a_1, a_2$ can be chosen in different ways (for example, both can be chosen to be very small or very large) giving rise to competing objectives.  Namely,
\begin{equation}  \label{objectives-competing}
\begin{array}{l}
x_5(t_f) := \ds\int\limits_0^{t_f} \Big(I(t) + L_2(t) + a_{11} \,u_1(t) + a_{12} \,u_2(t)\Big)\, dt\, , \\[3mm]
x_6(t_f) := \ds\int\limits_0^{t_f} \Big(I(t) + L_2(t) + a_{21} \,u_1(t) + a_{22} \,u_2(t)\Big)\, dt\,.
\end{array}
\end{equation}
with control weights $a_{11}, a_{12}, a_{21}, a_{22} \geq 0$, constitute two competing objective functionals.  Both functionals are given in Lagrange form. The standard method to obtain an optimal control problem of Bolza type is to introduce additional state variables $x_5$ and $x_6$ defined by
 \begin{equation}  \label{augmented_dyn}
\begin{array}{l}
\dot x_5(t) = I(t) + L_2(t) + a_{11}\, u_1(t) + a_{12}\, u_2(t)\,,\ \ x_5(0) = 0\,, \\[2mm]
 \dot x_6(t) = I(t) + L_2(t) + a_{21}\, u_1(t) + a_{22}\, u_2(t)\,,\ \ x_6(0) = 0\,.
 \end{array}
 \end{equation}
Denoting the (augmented) state vector by $x(t) =  (S(t),L_1(t),I(t),L_2(t),x_5(t),x_6(t)) \in \mathbb{R}^6$ and the control vector $u(t) := (u_1(t),u_2(t)) \in \mathbb{R}^2$,  the two competing objectives in the general problem (P) are given by
\[
\varphi_1 (x(t_f),t_f) = x_5(t_f) =: F_1(x,u)\quad\mbox{and}\quad
\varphi_2 (x(t_f),t_f) = x_6(t_f) =: F_2(x,u)\,,
\]
where $F_1(x,u)$ and $F_2(x,u)$ denote the two functionals in Lagrange form.

The bi-objective TB problem is now in the same form as Problem~(OCP) and, in particular, Problem~(OCPsd).  The decision maker's objective for this problem will be to minimize the distance to the origin of the value space.  We therefore choose
\[
\varphi_0(x^w,u^w,t_f^w) := \varphi_1^2(x^w(t_f),t_f^w) + \varphi_2^2(x^w(t_f),t_f^w)\,,
\]
Our aim is to solve Problem~(OPF), to determine a scalar $w\in(0,1)$ with $w_1 := w$ and $w_2 := 1-w$ that results in the best Pareto solution in the sense that $\varphi_0(\cdot,\cdot,\cdot)$ is minimized, subject to the solution of Problem~(OCP$_w$).

Next we focus on the solution of Problem~(OCP$_w$): We aim to find a pair of functions $(x,u)  \in W^{1,\infty}([0,t_f],\mathbb{R}^6) \times L^\infty([0,t_f],\mathbb{R}^2)$ that minimizes the parameter $\alpha$ subject to the time-delayed dynamics \eqref{dynamics} and the auxiliary dynamics \eqref{augmented_dyn}, initial conditions \eqref{initial-condition}, control constraints \eqref{control-constraint} and auxiliary weighted inequalities involving $\varphi_1$ and $\varphi_2$.

%%------------------------------------------------------------------
We consider the necessary optimality conditions for the 
time-delayed optimal control problem (OCP$_w$);
 see G\"ollmann and Maurer \cite{Goellmann-Maurer-2014}, Vinter \cite{Vinter}.
For this purpose we introduce the delayed state variable
 $y_3(t) = x_3(t-d_I) = I(t-d_I)$ and  delayed control variables $v_k(t) = u_k(t-d_{u_k})$, $k=1,2$.
 Denoting the adjoint variable vector by $\lambda(t) := (\lambda_S(t),\lambda_{L_1}(t),\lambda_I(t),\lambda_{L_2}(t),\lambda_5(t),\lambda_6(t)) \in \mathbb{R}^6$
 the Hamiltonian or Pontryagin function is given by
 %%-----------------------------------------------------------------
 \begin{equation}
 \label{Hamiltonian}
 \begin{array}{rl}
 H(x,y_3,\lambda,u_1,v_1,u_2,v_2) = &  
 \lambda_s\, (\mu N - \frac{\beta}{N} I S - \mu S) \\[2mm]
   &\ \ \, +\ \lambda_{L_1}\, (\frac{\beta}{N} I\left( S + \sigma L_2 + \sigma_R R\right) - \left(\delta + \tau_1
       + \epsilon_1 v_1 + \mu\right)L_1) \\[2mm]
&\ \ \, +\ \lambda_I\, (\,\phi \,\delta L_1 + \omega L_2 + \omega_R\, R - \tau_0 y_3 + \mu I ) \\[2mm]
&\ \ \, +\ \lambda_{L_2}\, ((1 - \phi) \delta L_1 - \sigma \frac{\beta}{N} I L_2
- (\omega + \epsilon_2 v_2 + \tau_2 + \mu)L_2 ) \\[2mm]
&\ \ \, +\ \lambda_5\,(I + L_2 + a_{11} u_1 + a_{12} u_2) \\[2mm]
&\ \ \, +\ \lambda_6\,(I + L_2 + a_{21} u_1 + a_{22} u_2)\,,
 \end{array}
 \end{equation}
 %%---------------------------------------------
where $R$ is given as in \eqref{R}.
 The Minimum Principle  \cite{Goellmann-Maurer-2014,Vinter} yields the adjoint equations
 \begin{equation*}
 \label{adjoint-equation}
 \dot{\lambda}_S(t) = -\frac{\partial H}{\partial S}[t] , \quad
 \dot{\lambda}_{L_1}(t) = -\frac{\partial H}{\partial L_1}[t] , \quad
 \dot{\lambda}_{L_2}(t) = -\frac{\partial H}{\partial L_2}[t] ,
 \end{equation*}
 \[
\dot{\lambda}_{x_5}(t) = -\frac{\partial H}{\partial x_5}[t] = 0\,, \quad
\dot{\lambda}_{x_6}(t) = -\frac{\partial H}{\partial x_6}[t] = 0\,,
 \]
 and the \textit{advanced} adjoint equation
  \begin{equation*}
 \label{adjoint-equation-I}
 \dot{\lambda}_I(t) = -\frac{\partial H}{\partial I}[t]
  - \chi_{[0,t_f-d_I]}(t) \frac{\partial H}{\partial I}[t+d_I]\,,
 \end{equation*}
where the argument $[t]$ stands for evaluating all arguments at time $t$.  We note that $\lambda_5^w(t) = \overline{\lambda}_5^w$ and $\lambda_6^w(t) = \overline{\lambda}_5^w$, constants, for any fixed $w\in[0,1]$.  In the last equation, the term $\chi_{[0,t_f-d_I]}(t)$ denotes the characteristic function of the interval $[0,t_f-d_I]$ at time $t$.
The minimization of the Hamiltonian with respect to the controls $u_1,u_2$ and delayed controls $v_1,v_2$ involves the {\em switching functions} $\sigma_k(t)$ for $k=1,2$:
\begin{equation}
\label{switching-functions}
\begin{array}{rcl}
\sigma_k^w(t) &=&  \ds\frac{\partial H}{\partial u_k}[t]
 +   \chi_{[0,t_f-d_{u_k}]}(t) \frac{\partial H}{\partial v_k}[t+d_{u_k}]
\\[4mm]
&=& \left\{
 \begin{array}{lcl}
  a_{1k}\overline{\lambda}_5^w + a_{2k}\overline{\lambda}_6^w 
  -\ \epsilon_k \lambda_{L_k}^w(t+d_{u_k})  L_k^w(t+d_{u_k})\,,
  & \mbox{if} & 0 \leq t \leq t_f-d_{u_k}\,,\\[1mm]
 a_{1k}\overline{\lambda}_5^w + a_{2k}\overline{\lambda}_6^w\,, & \mbox{if} &  t_f-d_{u_k} \leq t \leq t_f\,.
 \end{array}
\right .
\end{array}
\end{equation}
As in the Rayleigh problem, the superscript ``$w$'' above denotes dependence on the scalarization parameter/weight $w$.  Then the  controls minimizing the Hamiltonian are characterized by the switching conditions (control law)
\begin{equation}
\label{control-law}
u_k^w(t) = \left\{ \begin{array}{rcl}
          0\,, & \mbox{ if} & \sigma_k^w(t) > 0\,, \\
          1\,, & \mbox{ if} & \sigma_k^w(t) < 0\,,
          \end{array} \right . \qquad k=1,2.
\end{equation}
for all $w\in[0,1]$.  In particular, for positive weights $a_1 >0$, $a_2 >0$,  the switching functions \eqref{switching-functions}
and the control law \eqref{control-law} imply
$$
u_k^w(t) = 0 \quad \forall \; t_f-d_{u_k} \leq t \leq t_f \,,
$$
for all $w\in[0,1]$.

In what follows we choose the control weights as $a_{11}=a_{12} = 10$ (small) and $a_{21}=a_{22} = 1000$ (large) in the objective functionals $\varphi_1$ and $\varphi_2$.

\begin{figure}[t!]
\begin{center}
\includegraphics[width=120mm]{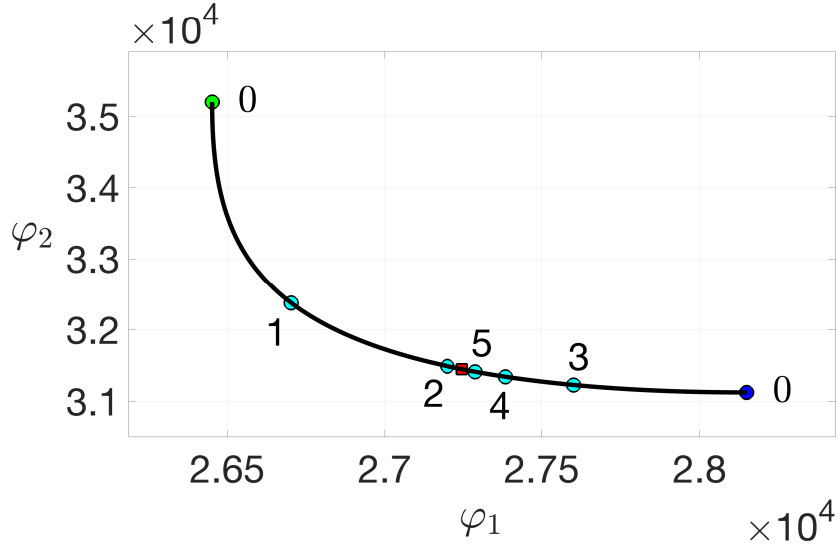}
\end{center}
\caption{\sf TB problem---Pareto front, and iterations of Algorithm~\ref{algo1}: Master solution is depicted by a (red) square and iterates by (light blue) circles..}
\label{fig:TB_Pareto}
\end{figure}

For Problem~(OCP$_w$) written for the TB problem, we have chosen the utopia vector as $(\beta_1^*, \beta_2^*) = (0, 0)$.  Figure~\ref{fig:TB_Pareto} depicts the Pareto front for the TB problem we consider here.  The plot also displays the iterations of Algorithm~\ref{algo1}.  The TB problem is discretized using the trapezoidal rule, the number of grid points is set to be $N = 5000$, and the Ipopt's tolerance to $10^{-10}$, so as to get solutions for $w$ accurate at least up to four decimal places~(dp).

The essential interval in this case is found to be $[w_0,w_f] = [0.5251, 0.5709]$, with
\[
(\varphi_1^{w_0},\varphi_2^{w_0}) = (28155, 31133)\quad\mbox{and}\quad
(\varphi_1^{w_f},\varphi_2^{w_f}) = (26459, 35205)\,,
\]
where $\varphi_i^{w} := \varphi_i(x^{w}(t_f),t_f^{w})$, $i = 1, 2$, with $w = w_0$ or $w_f$, or as will be the case below, $w = w^*$.
Optimization over the Pareto front results in $w^* = 0.5358$, after 10 iterations of Algorithm~\ref{algo1}, yielding 
\[
\varphi_0^{w^*} = 41621 \quad\mbox{and}\quad
(\varphi_1^{w^*},\varphi_2^{w^*}) = (27255,  31455)\,.
\]
In Figure~\ref{fig:TB_Pareto} only five iterations are displayed (labelled 1--5) for clarity in viewing.  The Pareto (master) solution with $w = w^*$ is represented by a square.

\begin{figure}[t!]
\begin{minipage}{75mm}
\begin{center}
\hspace*{-5mm}
\includegraphics[width=85mm]{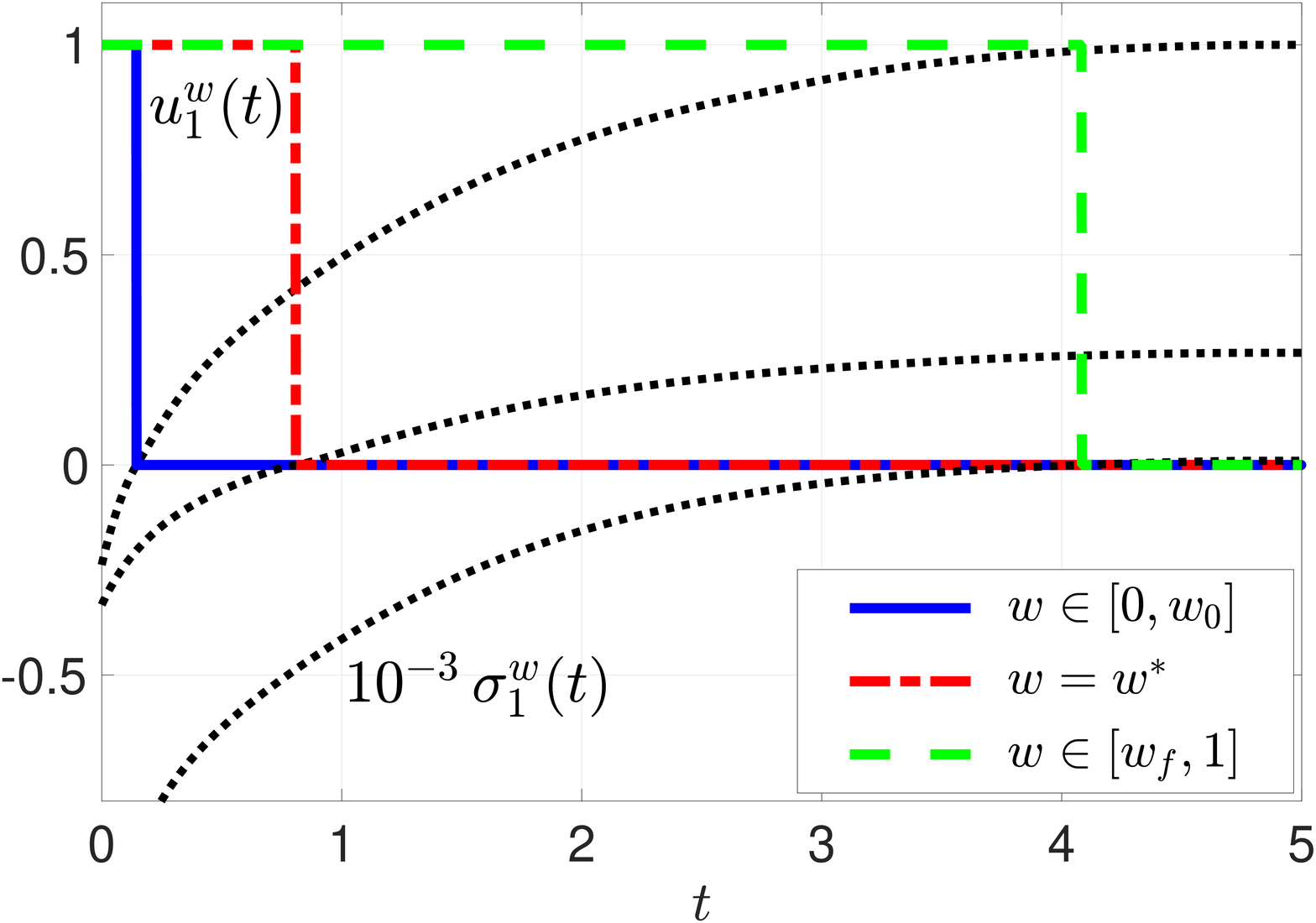} \\[3mm]
{\sf\small (a) Control variable $u_1^w$~\eqref{control-law} and scaled \\ switching function $\sigma_1^w$~\eqref{switching-functions} superposed.}
\end{center}
\end{minipage}
\begin{minipage}{85mm}
\begin{center}
\includegraphics[width=85mm]{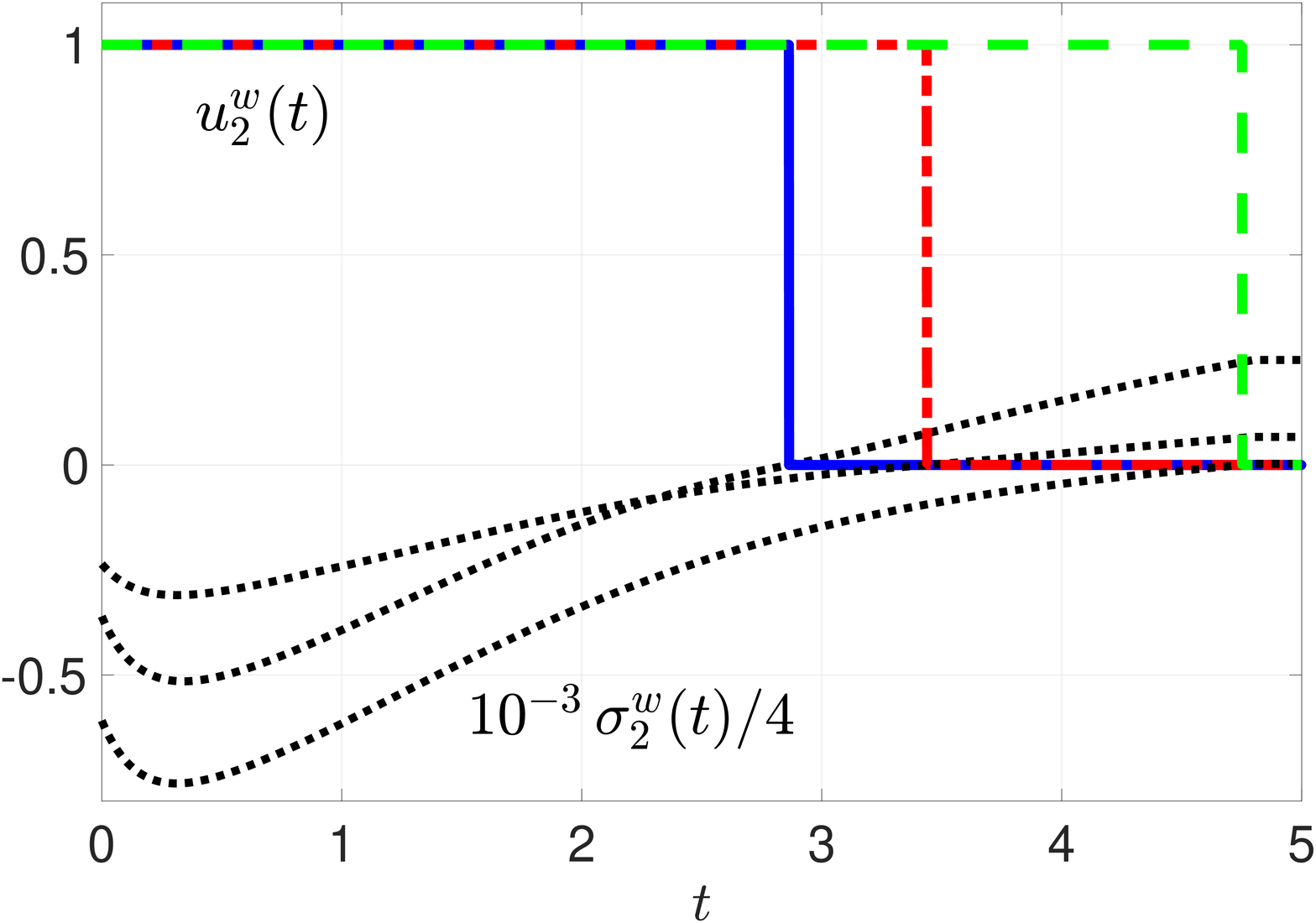} \\[3mm]
{\sf\small (b) Control variable $u_2^w$~\eqref{control-law} and scaled \\ switching function $\sigma_2^w$~\eqref{switching-functions} superposed.} 
\end{center}
\end{minipage}
\
\caption{\sf TB problem---Boundary Pareto solutions, corresponding to $w_0 = 0.5251$ and $w_f = 0.5709$, are shown with (blue) solid curves and (green) dashed curves, respectively.  Master Pareto solution, corresponding to $w^* = 0.5358$, is shown with dashed-and-dotted (red) curves.}
\label{fig:TB_controls}
\end{figure}

The numerical Pareto-optimal control variable solutions $u_1^w(\cdot)$ and $u_2^w(\cdot)$ are presented in~Figures~\ref{fig:TB_controls}(a)--(b) for $w = w_0, w^*, w_f$.  As with Rayleigh, one of the boundary Pareto-optimal solutions is shown using solid (blue) curves for $w = w_0$, the same solution for all $w\in[0,w_0]$.  The other boundary Pareto-optimal solution for $w = w_f$, which holds for all $w\in[w_f,1]$, is shown using dashed (green) curves.  Both of the control solutions are of bang--bang type (as required by \eqref{control-law}), with one switching (the number of switchings not dictated by \eqref{control-law} alone).  The master Pareto solution is given for $w = w^*$ using dashed-and-dotted (red) curves, in which the controls are also of bang--bang type with one switching.

The switching functions for each control and case, $\sigma_k^w(\cdot)$, $k = 1,2$, scaled as indicated, are plotted with (black) dotted curves and superposed with the control plots in Figures~\ref{fig:TB_controls}(a)--(b).  We remind that, by using~\eqref{switching-functions} (recall that discrete approximations of $\lambda^w_{L_k}(t)$, $k = 1,2$, $\lambda^w_5(t)$ and $\lambda^w_6(t)$ can readily be obtained as constraint multipliers from AMPL), one verifies the optimality condition in~\eqref{control-law}.

In each strategy,  the two control efforts are ``on'' until the times $t_{s_k}^w$, $k = 1,2$, at which the respective $u_k^w(\cdot)$ is switched ``off'' (down to zero).  These types of bang--bang controls are also referred to as {\em on--off controls}.  In Table~\ref{tab:TB} the switching times for the boundary as well as the optimal weights are listed.  Under these controls, the resulting terminal values of the state variables are also listed in Table~\ref{tab:TB}.  The plots of these variables are not provided as they are difficult to distinguish at earlier times (as expected) and that they become distinguishable/comparable only near the terminal time.

Under the controls minimizing $x_5(t_f)$ (with $w = w_f = 0.5709$ and minimum $x_5^{w_f}(t_f) = 26459$) the number of persistent latent individuals $L_2(t_f)$ turns out to be about 419 (in a population of 30000).  This number is more than doubled to 864 if $x_6(t_f)$ is minimized (with $w = w_0 = 0.5709$ and minimum $x_6^{w_0}(t_f) = 31133$).  The optimal Pareto solution minimizing the distance in value space to the origin yields with $w = w^* = 0.5358$ the optimal $L_2(t_f)$ as 748.

 \begin{table}[t]
\centering
{\small
\begin{tabular}{cccccrrrrr}
 Scalarization & \multicolumn{2}{|c|}{Functional values} 
 & \multicolumn{2}{|c|}{Switching times}
 & \multicolumn{5}{c|}{Terminal state values} \\  
 \cline{2-3}  \cline{4-5}  \cline{6-10}
weight $w$ 
& \multicolumn{1}{|c}{$x_5^w(t_f)$} & \multicolumn{1}{c}{$x_6^w(t_f)$} 
& \multicolumn{1}{|c}{$t_{s_1}^w$} & \multicolumn{1}{c}{$t_{s_2}^w$} 
& \multicolumn{1}{|c}{$S^w(t_f)$} & \multicolumn{1}{c}{$L_1^w(t_f)$}  
& \multicolumn{1}{c}{$I^w(t_f)$} & \multicolumn{1}{c}{$L_2^w(t_f)$}  
& \multicolumn{1}{c|}{$R^w(t_f)$}  \\ \hline\hline
$w_0 = 0.5251$\,:
 & 28155 & 31133 & 0.145 & 2.864 & 1193.1 & 28.2 & 13.3 & 864.0 & 27901.4 \\
$w^* = 0.5358$\,: & 27255 & 31455 & 0.809 & 3.439 & 1205.8 & 27.5 & 13.0 & 747.6 & 28006.1 \\
$w_f = 0.5709$\,: & 26459 & 35205 & 4.083 & 4.752 & 1238.2 & 23.8 & 11.2 & 419.3 & 28307.5 \\
     \hline
    \end{tabular}}
\caption{\sf\small TB problem.}
\label{tab:TB}
\end{table}

\section{Conclusion}

We have proposed an algorithm to solve the problem of optimization over the Pareto front.  The algorithm employs bisection method which starts with an essential interval of weights of the Chebyshev scalarization.  It is applicable to a wide range of optimal control problems, including state- and control-constrained problems with time delay.  Numerical solution of two challenging optimal control problems has demonstrated the effectiveness of the algorithm.

The main motive behind the algorithm we have proposed is that one can find the optimal solution minimizing a master objective functional without having to construct the Pareto front.  The algorithm solves the challenging optimal control problem (OCP$_w$) a relatively smaller number of times than the case  of constructing the Pareto front.  In the examples we have studied the algorithm had to solve (OCP$_w$) 20 to 30 times.  On the other hand, without the algorithm we propose, it is necessary to construct the Pareto front by solving (OCP$_w$) thousands of times in order to obtain the same solution with the same computational accuracy.

The proposed algorithm can be improved/modified in various ways.  For example, scalarization techniques other than Chebyshev might be employed; see for example~\cite{BurKayRiz2013, BurKayRiz2017} and the references therein.  Bisection method might be replaced by methods with higher convergence rates, for example {\em regula falsi} and secant methods (see \cite{BurFai2011}), at the expense of approximating higher order derivatives of course, although the latter would make the algorithm applicable to problems with more than just two objective functionals.

%\newpage

\end{document}